\documentclass[11pt,a4paper]{article}
\usepackage{amsmath}
\usepackage{amsfonts}
\usepackage{amsthm}
\usepackage{graphicx}
\usepackage{geometry}
\usepackage{multirow}
\usepackage{algorithmic}
\usepackage{algorithm}
\usepackage{hhline}
\numberwithin{algorithm}{section}

\begin{document}

\markboth{Zhiwei Qin, Donald Goldfarb, and Shiqian Ma}{Optimization Methods and Software}


\title{An Alternating Direction Method for Total Variation Denoising}

\author{Zhiwei (Tony) Qin \footnote{Department of Industrial Engineering and Operations Research, Columbia University, New York, NY 10027.  zq2107@columbia.edu.} \and Donald Goldfarb \footnote{Department of Industrial Engineering and Operations Research, Columbia University, New York, NY 10027.  goldfarb@columbia.edu.} \and Shiqian Ma \footnote{Department of Systems Engineering and Engineering Management, The Chinese University of Hong Kong, Hong Kong, China.  sqma@se.cuhk.edu.hk.}}

\maketitle

\newtheorem{alg}{Algorithm}[section]
\newtheorem{thm}{Theorem}[section]
\newtheorem{cor}{Corollary}[section]
\newtheorem{lem}{Lemma}[section]
\newtheorem{rem}{Remark}[section]
\newcommand{\xsupk}[2]{#1^{(#2)}}
\newcommand{\usupk}{\xsupk{u}{k}}
\newcommand{\xkpone}{\xsupk{x}{k+1}}
\newcommand{\ykpone}{\xsupk{y}{k+1}}
\newcommand{\gammakpone}{\xsupk{\gamma}{k+1}}
\newcommand{\gammak}{\xsupk{\gamma}{k}}
\newcommand{\bbR}{\mathbb{R}}

\begin{abstract}
We consider the image denoising problem using total variation (TV) regularization.  This problem can be computationally challenging to solve due to the non-differentiability and non-linearity of the regularization term.  We propose an alternating direction augmented Lagrangian (ADAL) method, based on a new variable splitting approach that results in subproblems that can be solved efficiently and exactly.  The global convergence of the new algorithm is established for the anisotropic TV model.  For the isotropic TV model, by doing further variable splitting, we are able to derive an ADAL method that is globally convergent. We compare our methods with the split Bregman method \cite{goldstein2009split},which is closely related to it, and demonstrate their competitiveness in computational performance on a set of standard test images.
\bigskip


\end{abstract}

\section{Introduction}
In signal processing, total variation (TV) regularization is a very popular and effective approach for noise reduction and has a wide array of applications in digital imaging.  The total variation is the integral of the absolute gradient of the signal.  Using TV regularization to remove noise from signals was originally proposed in \cite{rudin1992nonlinear} and is based on the observation that noisy signals have high total variation.  By reducing the total variation of a noisy signal while keeping the resulting signal close to the original one removes noise while preserving important details such as sharp edges.  Other existing denoising techniques include median filtering and Tikhonov-like regularization, $\|u\|_{TIK} := \sum_i (\nabla_xu)_i^2 + (\nabla_yu)_i^2$ of the desired solution $u$, where $\nabla_x$ and $\nabla_y$ are defined below in \eqref{eq:diff_u}.  It is known that they tend to smooth away important texture details along with the noise \cite{strong2003edge, wang2008new}.

For a 2-D signal $u \in \mathbb{R}^{n\times m}$, such as an image, the total variation $\|u\|_{TV}$ \cite{rudin1992nonlinear} of $u$ can be defined anisotropically or isotropically:
\begin{equation}\label{tv_norms}
    \|u\|_{TV} := \left\{
                   \begin{array}{ll}
                     \sum_i |(\nabla_x u)_i| + |(\nabla_y u)_i|, & \hbox{(Anisotropic);} \\
                     \sum_i \sqrt{(\nabla_xu)_i^2 + (\nabla_yu)_i^2}, & \hbox{(Isotropic).}
                   \end{array}
                 \right.
\end{equation}
Concisely, $\|u\|_{TV}$ can be expressed as $\sum_{i=1}^{nm}\|D_i u\|_p$, where $D_i u \in \mathbb{R}^2$ denotes the discrete gradient of $u$ at pixel $i$.  Hence, $\|u\|_{TV}$ is isotropic when $p=2$ and is anisotropic when $p=1$.
TV denoising (also called ROF (Rudin-Osher-Fatemi) denoising) corresponds to solving the following optimization problem,
\begin{equation}\label{eq:tvdn}
    \min_u \lambda\sum_{i=1}^{nm}\|D_i u\|_p + \frac{1}{2}\|u-b\|^2,
\end{equation}
where $p = 1$ or 2; $b \in \mathbb{R}^{n\times m}$ is the noisy image, and the solution $u$ is the desired denoised image.  $\|\cdot\|$ without a subscript denotes the $l_2$-norm.
We assume that all 2-D images are in column-major vectorized form; hence, if one-dimensional index of $(i,j)$ is $k$ and $1 \leq i \leq n, 1 \leq j \leq m$, the elements of $\nabla u$ are given by
\begin{equation}\label{eq:diff_u}
    D_k u \equiv  \left(
                 \begin{array}{c}
                   u_{k+1}-u_k \\
                   u_{k+n}-u_k \\
                 \end{array}
               \right) = [\nabla u]_{ij} = \left(
                           \begin{array}{c}
                             \nabla_x u \\
                             \nabla_y u \\
                           \end{array}
                         \right)_{ij}.
\end{equation}
The anisotropic TV model that we consider in this paper is the four-neighbor form.  Algorithms for anisotropic TV denoising for other different sets of neighbors are presented in \cite{goldfarb2009parametric}.

Due to the non-differentiability and non-linearity of the TV term in problem \eqref{eq:tvdn}, this problem can be computationally challenging to solve despite its simple form.  Hence, much effort has been devoted to devise efficient algorithms for solving it.  A number of references are provided in Section 1 of \cite{goldstein2009split}.  In addition, Chambolle's algorithm \cite{chambolle2004algorithm} solves problem \eqref{eq:tvdn} with the isotropic TV-norm.

The approach that we develop in this paper for solving problem \eqref{eq:tvdn} is based on variable splitting followed by the application, to the resulting constrained minimization problem, of an alternating minimization algorithm (specifically, in our case, the alternating direction augmented Lagrangian (ADAL) method).  In contrast with previously proposed variable splitting approaches for solving problem \eqref{eq:tvdn}, our approach introduces two sets of auxiliary variables, one to replace the solution image $u$ and one to replace the vector of gradients ($D_1 u, \cdots, D_{nm} u$).  When the ADAL method is applied to the constrained optimization problem that is derived from this variable splitting, the resulting subproblems that must be solved at each iteration can be solved easily and exactly.  Moreover, for the anisotropic TV version of problem \eqref{eq:tvdn}, convergence of our algorithm can be proved, and for both the anisotropic and isotropic TV models, preliminary numerical experiments indicate that the number of iterations required to obtain an accurate solution is quite small.  By introducing a third set of auxiliary variables, we are also able to derive an ADAL method for the isotropic TV model with guaranteed convergence.

Before outlining the organization of the remaining sections of the paper that contain our main results, let us first review three previously proposed methods that use split-variable alternating minimization approaches.  All of these methods apply to the slightly more general TV-based denoising/deblurring problem
\begin{equation}\label{eq:tv_deblur}
    \min_{u} \lambda\sum_i \|D_i u\|_p + \frac{1}{2}\|Ku - b\|^2,
\end{equation}
where $p$ is either 1 or 2, and $K$ is a blurring (or convolution) operator.

\subsection{Closely Related Methods}\label{sec:related_work}
In the straightforward variable splitting approach proposed in \cite{afonso2010fast}, a vector of auxiliary variables $w$ is introduced to replace $u$ in the non-differentiable TV term in \eqref{eq:tv_deblur}:
\begin{eqnarray}\label{eq:tv_deblur_split}
  \min_{u,w} && \lambda\sum_i\|D_i w\|_p + \frac{1}{2}\|Ku-b\|^2 \\
  \nonumber s.t. && w = u.
\end{eqnarray}
The algorithm SALSA (Split-Augmented Lagrangian Shrinkage Algorithm) in \cite{afonso2010fast} then obtains a solution to problem \eqref{eq:tv_deblur} by applying the ADAL method to problem \eqref{eq:tv_deblur_split}, in which the non-differentiable TV term $\lambda\|\Phi(w)\|_p$, where $\Phi_i(w)\equiv D_i w$, has been decoupled from the quadratic fidelity term $R(u)\equiv \frac{1}{2}\|Ku-b\|^2$ in the objective function.  For the case of isotropic TV regularization, SALSA uses five iterations of Chambolle's algorithm to compute the corresponding Moreau proximal mapping.

In \cite{wang2008new}, variable-splitting combined with a penalty function approach is applied to problem \eqref{eq:tv_deblur} by introducing an auxiliary variable $d_i = D_i u \in \mathbb{R}^2$ for each pixel, yielding the following approximation to problem \eqref{eq:tv_deblur}
\begin{equation}\label{eq:quad_penalty}
    \min_{d,u} \lambda\sum_i \|d_i\|_1 + \frac{1}{2}\|Ku - b\|^2 + \frac{1}{2\mu}\sum_i \|d_i-D_i u\|^2.
\end{equation}
Problem \eqref{eq:quad_penalty} is then minimized alternatingly with respect to $w$ and $u$, with a continuation scheme that drives the penalty parameter $\frac{1}{\mu}$ gradually to a sufficiently large number.  This method is extended in \cite{yang2009fast, yang2009efficient} to solve the multi-channel (color) image deblurring problem.  In \cite{yang2009efficient}, the TV regularization with 1-norm fidelity (TVL1) model
\begin{equation*}
    \min_u \lambda\sum_i \|D_i u\|_p + \|Ku - b\|_1
\end{equation*}
is considered.  The same approach has also been applied to reconstruct signals from partial Fourier data in the compressed sensing context \cite{yang2010fast}.  These methods take full advantage of the structures of the convolution operator and the finite difference operator so that the subproblems can be solved exactly and efficiently, which is important for fast convergence.  A downside to this quadratic penalty approach is that when $\frac{1}{\mu}$ is very large, problem \eqref{eq:quad_penalty} becomes ill-conditioned and numerical stability becomes an issue.

Although our algorithm is closely related to the algorithms in \cite{afonso2010fast} and \cite{wang2008new} described above, it is even more closely related to the split Bregman method \cite{goldstein2009split}, which is an application of the variable splitting approach in \cite{wang2008new} to the Bregman iterative regularization method \cite{osher2006iterative}.  The Bregman iterative regularization method was first introduced in \cite{osher2006iterative} as a better (iterative) approach to the TV denoising/deblurring problem \eqref{eq:tv_deblur} than directly applying an iterative solver to it.  Subsequently, this method was extended in \cite{strong2003edge} to the solution of $l_1$-minimization problems that arise in compressed sensing and in \cite{ma2009fixed} to nuclear norm minimization problems that are convex relaxations of matrix rank minimization problems.

The Bregman distance associated with a convex function $E(\cdot)$ between $u$ and $v$ is defined as
\begin{equation*}
    D^p_E(u,v) := E(u) - E(v) - p^T(u-v),
\end{equation*}
where $p \in \partial E(v)$ and $\partial E(v)$ denotes the subdifferential of $E(\cdot)$ at the point $v$.
The Bregman iteration for the unconstrained minimization problem
\begin{equation*}
    \min_u E(u) + \frac{1}{\mu} H(u),
\end{equation*}
where both functions $E(\cdot)$ and $H(\cdot)$ are convex,
is
\begin{eqnarray}
  \nonumber\xsupk{u}{k+1} &=& \arg\min_u D^p_E(u,\xsupk{u}{k}) + \frac{1}{\mu} H(u) \\
   &=& \arg\min_u E(u) - (u-\xsupk{u}{k})^T\xsupk{p}{k} + \frac{1}{\mu} H(u), \label{eq:bregman_iter1}\\
  \xsupk{p}{k+1} &=& \xsupk{p}{k} - \nabla H(\xsupk{u}{k+1}). \label{eq:bregman_iter2}
\end{eqnarray}
Superscripts denote iteration indices to differentiate between the values of the variables for the current iteration from those computed at the previous iteration.
With the introduction of an auxiliary variable $d$ as in \cite{wang2008new}, the TV denoising/deblurring problem \eqref{eq:tv_deblur} can be reformulated as the constrained problem
\begin{eqnarray}\label{eq:tvdn_sb_constr}
  \min_{u,d} && \lambda\|d\|_1 + R(u) \\
  \nonumber s.t. && d = \Phi(u),
\end{eqnarray}
where
$R(u) = \frac{1}{2}\|Ku-b\|^2$, and $\Phi(u) = \left(
                                      \begin{array}{c}
                                        \nabla_x u \\
                                        \nabla_y u \\
                                      \end{array}
                                    \right)
$.  Now, converting problem \eqref{eq:tvdn_sb_constr} into an unconstrained problem (by penalizing $\|d-\Phi(u)\|^2$), we obtain
\begin{equation*}
    \min_{u,d} \lambda\|d\|_1 + R(u) + \frac{1}{2\mu}\|d-\Phi(u)\|^2,
\end{equation*}
(this is the same as problem \eqref{eq:tv_deblur_split}). Then, applying the general Bregman iteration \eqref{eq:bregman_iter1}-\eqref{eq:bregman_iter2} with $E(u,d) = \lambda\|d\|_1 + R(u)$ and $H(u,d) = \|d-\Phi(u)\|^2$, we obtain after simplification the following specific Bregman iteration:
\begin{eqnarray}
  (\xsupk{u}{k+1}, \xsupk{d}{k+1}) &=& \min_{u,d} \lambda\|d\|_1 + R(u) + \frac{1}{2\mu}\|d-\Phi(u)-\xsupk{r}{k}\|^2, \label{eq:sb_subprob} \\
  \xsupk{r}{k+1} &=& \xsupk{r}{k} + (\Phi(\xsupk{u}{k+1}) - \xsupk{d}{k+1}), \label{eq:bregman_update}
\end{eqnarray}
with $\xsupk{r}{0} = 0$.
In \cite{goldstein2009split}, an approximate solution to \eqref{eq:sb_subprob} was proposed by alternatingly minimizing the right-hand-side of \eqref{eq:sb_subprob} with respect to $u$ and $d$ once.  This yields the following Split Bregman algorithm (Algorithm \ref{alg:splitbregman}).  \footnote{The Split Bregman method in its original form in \cite{goldstein2009split} has an inner loop.  We consider the simplified form, which was used to solve TV denoising problems in \cite{goldstein2009split}.} For notational conciseness, the superscripts are suppressed in the main steps.
\begin{algorithm}
\caption{SplitBregman}
\begin{algorithmic}[1]\label{alg:splitbregman}
\STATE Given $\xsupk{u}{0}$, $\xsupk{d}{0}$, and $\xsupk{r}{0}$.
\FOR{$k = 0,1,\cdots,K$}
    \STATE $u \gets \min_u R(u) + \frac{1}{2\mu}\|d- \Phi(u) - r\|^2$ \label{line:sb_u}
    \STATE $d \gets \min_d \lambda\|d\|_1 + \frac{1}{2\mu}\|d - \Phi(u) - r\|^2$ \label{line:sb_d}
    \STATE $r \gets r + (\Phi(u) - d)$
\ENDFOR
\RETURN $u$
\end{algorithmic}
\end{algorithm}
As is well known (e.g., see \cite{yang2009alternating}), the Bregman iterative algorithm \eqref{eq:sb_subprob}-\eqref{eq:bregman_update} is equivalent to applying the augmented Lagrangian method \cite{hestenes1969multiplier, powell1972nonlinear} to solve problem \eqref{eq:tvdn_sb_constr}.  Hence, the split-Bregman algorithm is equivalent to applying the ADAL method to \eqref{eq:tvdn_sb_constr}.

\subsection{Organization of The Paper}
The outline of the rest of the paper is as follows.  We first briefly review the ADAL method and its applications to linearly constrained optimization problems that arise from variable splitting in Section \ref{sec:adal}.  In Section \ref{sec:adal_tvdn}, we describe our proposed variable-splitting alternating direction augmented Lagrangian method for the anisotropic TV-model and prove its global convergence in Section \ref{sec:conv_aniso}.  We then discuss in Sections \ref{sec:iso} and \ref{sec:contrast_sb_tvdn} the isotropic case and the difference between our algorithm and the split Bregman method, respectively.  In Section \ref{sec:iso_conv}, we present a globally convergent variable-splitting ADAL variant for the isotropic TV-model. In Section \ref{sec:exp}, we compare our algorithms against the split Bregman method on a set of standard test images and demonstrate the effectiveness of our methods in terms of denoising speed and quality.

\section{The Alternating Direction Augmented Lagrangian Method}\label{sec:adal}

The ADAL method is also known as the alternating direction method of multipliers (ADMM) and was first proposed in the 1970s \cite{gabay1976dual, glowinski1975adal}.  It belongs to the family of the classical augmented Lagrangian (AL) method \cite{powell1972nonlinear, rockafellar1973multiplier, hestenes1969multiplier}, which iteratively solves the linearly constrained problem
\begin{eqnarray}\label{eq:lin_constr_prob}
  \min_x && F(x) \\
  \nonumber s.t. && Ax = b.
\end{eqnarray}
The augmented Lagrangian of problem \eqref{eq:lin_constr_prob} is $\mathcal{L}(x,\gamma) = F(x) + \gamma^T(b-Ax) + \frac{1}{2\mu}\|Ax-b\|^2$, where $\gamma$ is the Lagrange multiplier and $\mu$ is the penalty parameter for the quadratic infeasibility term.  The AL method minimizes $\mathcal{L}(x,\gamma)$ followed by an update to $\gamma$ in each iteration as stated in the following algorithm.  We denote by $k_{max}$ the user-defined maximum number of iterations or the number of iterations required to satisfy the termination criteria.
\begin{algorithm}
\caption{AL (Augmented Lagrangian method)}
\begin{algorithmic}[1]\label{alg:al}
\STATE Choose $\xsupk{\gamma}{0}$.
\FOR{$k = 0,1,\cdots,k_{max}$}
    \STATE $x \gets$ $\arg\min_x\mathcal{L}(x, \gamma)$ \label{line:solve_auglag}
    \STATE $\gamma \gets \gamma - \frac{1}{\mu}(Ax - b)$
\ENDFOR
\RETURN $x$
\end{algorithmic}
\end{algorithm}

For a structured unconstrained problem
\begin{equation}\label{eq:struct_uncon_prob}
    \min_x F(x) \equiv f(x) + g(Ax),
\end{equation}
where both functions $f(\cdot)$ and $g(\cdot)$ are convex, we can decouple the two functions by introducing an auxiliary variable $y$ and transform problem \eqref{eq:struct_uncon_prob} into an equivalent linearly constrained problem
\begin{eqnarray}\label{eq:struct_lin_con}
  \min_{x,y} && f(x) + g(y) \\
  \nonumber s.t. && Ax = By,
\end{eqnarray}
where for \eqref{eq:struct_uncon_prob}, $B=I$.  Henceforth, we consider the more general case of problem \eqref{eq:struct_lin_con}.
The augmented Lagrangian function for problem \eqref{eq:struct_lin_con} is
\begin{equation*}
    \mathcal{L}(x,y,\gamma) = f(x) + g(y) + \gamma^T(By-Ax) + \frac{1}{2\mu}\|By-Ax\|^2.
\end{equation*}
Exact joint minimization of $\mathcal{L}(x,y,\gamma)$ with respect to both $x$ and $y$ is usually difficult.  Hence, in practice, an inexact version of the AL method (IAL) is often used, where $\mathcal{L}(x,y,\gamma)$ is minimized only approximately.  Convergence is still guaranteed in this case, as long as the subproblems are solved with increasing accuracy \cite{rockafellar1973multiplier}.

ADAL (Algorithm \ref{alg:adal} below) is a particular case of IAL in that it finds the approximate minimizer of $\mathcal{L}(x,y,\gamma)$ by alternatingly optimizing with respect to $x$ and $y$ once.  This is often desirable because joint minimization of $\mathcal{L}(x,y,\gamma)$ even approximately can be hard.
\begin{algorithm}
\caption{ADAL (ADMM)}
\begin{algorithmic}[1]\label{alg:adal}
\STATE Choose $\xsupk{\gamma}{0}$.
\FOR{$k = 0,1,\cdots,k_{max}$}
    \STATE $x \gets \arg\min_{x}\mathcal{L}(x, y, \gamma )$
    \STATE $y \gets \arg\min_{y}\mathcal{L}(x, y, \gamma )$
    \STATE $\gamma \gets \gamma + \frac{1}{\mu}(By - Ax)$
\ENDFOR
\RETURN $x$
\end{algorithmic}
\end{algorithm}

The convergence of ADAL has been established for the case of two-way splitting as above.  This result, which is a modest extension of results in \cite{eckstein1992douglas}, is given in \cite{esser2009applications} and contained in the following theorem.
\begin{thm}\label{thm:admm_conv}
Consider problem \eqref{eq:struct_lin_con}, where both $f$ and $g$ are proper, closed, convex functions, and $A \in \mathbb{R}^{n\times m}$ and $B \in \mathbb{R}^{n\times l}$ have full column rank.  Then, starting with an arbitrary $\mu > 0$ and $x^0 \in \bbR^m, y^0 \in \mathbb{R}^l$, the sequence $\{x^k,y^k,\gamma^k\}$ generated by Algorithm \ref{alg:adal} converges to a primal-dual optimal solution pair $\big((x^*,y^*),\gamma^*\big)$ to problem \eqref{eq:struct_lin_con}, if \eqref{eq:struct_lin_con} has one.  If \eqref{eq:struct_lin_con} does not have an optimal solution, then at least one of the sequences $\{(x^k,y^k)\}$ and $\{\gamma^k\}$ diverges.
\end{thm}
It is known that $\mu$ does not have to decrease to a very small value (it can simply stay constant) in order for the method to converge to the optimal solution of problem \eqref{eq:struct_lin_con} \cite{nocedal, bertsekas1999nonlinear}.  Inexact versions of ADAL, where one or both of the subproblems are solved approximately have also been proposed and analyzed \cite{eckstein1992douglas, he2002new,yang2009alternating}.

The versatility and simple form of ADAL have attracted much attention from a wide array of research fields.  ADAL has been applied to solve group sparse optimization problems in \cite{deng2011group}, semidefinite programming problems in \cite{wen2010alternating} and matrix completion problems with nonnegative factors in \cite{xu2012alternating}. In signal processing/reconstruction, ADAL has been applied to sparse and low-rank recovery, where nuclear norm minimization is involved \cite{lin2010augmented, yuan2009sparse, shen2011augmented}, and to the $l_1$-regularized problems in compressed sensing \cite{yang2009alternating}.  ADAL-based algorithms have also been proposed to solve a number of image processing tasks, such as image inpainting and deblurring (SALSA and C-SALSA) \cite{afonso2010fast, afonso2009augmented, almeida2012deconvolving,tao2009alternating}, motion segmentation and reconstruction \cite{ramani2012splitting,zappella11simul}, in addition to denoising \cite{afonso2010fast,goldstein2009split,figueiredo2010restoration,steidl2010removing}.  In machine learning, ADAL and IAL-based methods have been successfully applied to structured-sparsity estimation problems \cite{qin2011structured} as well as many others \cite{boyd2010distributed}.

\section{Our Proposed Method}
\subsection{Application to Anisotropic TV Denoising}\label{sec:adal_tvdn}
We consider the anisotropic TV denoising model \eqref{eq:tvdn}.  The isotropic TV model will be considered in Section \ref{sec:iso}.  As in \cite{goldstein2009split}, we introduce auxiliary variables $d_x$ and $d_y$ for the discretized gradient components  $\nabla_x u$ and $\nabla_y u$ respectively.  Under reflective Neumann boundary conditions, $\nabla_x u = Du$, where the discretization matrix $D$ is an $(nm-m)\times nm$ block diagonal matrix, each of whose $m$ diagonal $(n-1)\times n$ rectangular blocks is upper bidiagonal with -1's on its diagonal and 1's on its super-diagonal.  For simplicity, henceforth we will assume that $n=m$.  Consequently, $\nabla_y u = Dv$, where $v = Pu$, and $P$ is a permutation matrix so that $v$ is the row-major vectorized form of the 2-D image. (Recall that $u$ is in the column-major form.)  Hence, problem \eqref{eq:tvdn} is equivalent to the following constrained problem
\begin{eqnarray}\label{eq:aniso_constr}
  \min_{d_x,d_y,u,v} && \lambda(\|d_x\|_1 + \|d_y\|_1) + \frac{1}{2}\|u-b\|^2 \\
  \nonumber s.t. && d_x = Du, \\
  \nonumber && d_y = Dv, \\
  \nonumber && v = Pu.
\end{eqnarray}
The augmented Lagrangian for problem \eqref{eq:aniso_constr} is
\begin{multline}\label{eq:aniso_al}
    \mathcal{L}(d_x,d_y,u,v,\mu) \equiv \frac{1}{2}\|u-b\|^2 + \lambda(\|d_x\|_1 + \|d_y\|_1) + \gamma_x^T(Du - d_x) + \gamma_y^T(Dv - d_y) + \gamma_z^T(Pu - v) \\
    + \frac{1}{2\mu_1}(\|Du-d_x\|^2 + \|Dv-d_y\|^2) + \frac{1}{2\mu_2}\|Pu-v\|^2.
\end{multline}
To minimize $\mathcal{L}$ with respect to $d = \left(
                                                 \begin{array}{c}
                                                   d_x \\
                                                   d_y \\
                                                 \end{array}
                                               \right)
$, we solve the subproblem
\begin{equation}\label{eq:aniso_sub_d}
    \min_{d_x,d_y} \lambda(\|d_x\|_1 + \|d_y\|_1) + \gamma_x^T(Du - d_x) + \gamma_y^T(Dv - d_y)
    + \frac{1}{2\mu_1}(\|Du-d_x\|^2 + \|Dv-d_y\|^2).
\end{equation}
Problem \eqref{eq:aniso_sub_d} is strictly convex and decomposable with respect to $d_x$ and $d_y$, so
the unique minimizer can be computed through two independent soft-thresholding operations
\begin{eqnarray*}
  d_x^* &=& \mathcal{T}(Du + \mu_1\gamma_x, \lambda\mu_1), \\
  d_y^* &=& \mathcal{T}(Dv + \mu_1\gamma_y, \lambda\mu_1),
\end{eqnarray*}
where the soft-thresholding operator $\mathcal{T}$ is defined componentwise as
\begin{equation*}
    \mathcal{T}(x,\lambda)_i := \max\{ |x_i| - \lambda, 0 \}\textrm{sign}(x_i).
\end{equation*}

To minimize $\mathcal{L}$ over $u$, we solve
\begin{equation}\label{eq:aniso_sub_u}
    \min_u \frac{1}{2}\|u-b\|^2 + \gamma_x^TDu + \frac{1}{2\mu_1}\|Du-d_x\|^2 + \gamma_z^TPu + \frac{1}{2\mu_2}\|Pu-v\|^2,
\end{equation}
which simplifies to the linear system
\begin{equation}\label{eq:aniso_sys_u}
    \left(D^TD + \left(\frac{\mu_1}{\mu_2}+\mu_1\right)I\right)u = \mu_1b + D^T(d_x - \mu_1\gamma_x) + P^T\left(\frac{\mu_1}{\mu_2}v - \mu_1\gamma_z\right).
\end{equation}
It is easy to verify that, since $\mu_1$ and $\mu_2$ are both positive scalars, the matrix on the left-hand-side of the above system is positive definite and tridiagonal.  Hence, \eqref{eq:aniso_sys_u} can be solved efficiently by the Thomas algorithm in $8nm$ flops \cite{golub1996matrix}.  We denote the solution to the above tridiagonal system by $u(d_x,v,\gamma_x,\gamma_z)$.

Similarly, the sub-problem with respect to $v$ simplifies to the tridiagonal system
\begin{equation}\label{eq:aniso_sys_v}
    \left( D^TD + \frac{\mu_1}{\mu_2}I \right)v = D^T(d_y - \mu_1\gamma_y) + \mu_2\gamma_z + \frac{\mu_1}{\mu_2}Pu.
\end{equation}
Its solution is denoted by $v(d_y,v,\gamma_y,\gamma_z)$.

With all the ingredients of the algorithm explained, we formally state this ADAL method in Algorithm \ref{alg:adal_tvdn} below.  Note that in line \ref{line:LM-update}, the vectors of Lagrange multipliers and scaled infeasibilities are combined into the vectors
\begin{equation*}
\gamma \equiv \left(
                                                                                                                              \begin{array}{c}
                                                                                                                                \gamma_x \\
                                                                                                                                \gamma_y \\
                                                                                                                                \gamma_z \\
                                                                                                                              \end{array}
                                                                                                                            \right) \; \textrm{and}\; \Delta \equiv \left(
                                                                                                                              \begin{array}{c}
                                                                                                                                \frac{1}{\mu_1}(Du - d_x) \\
                                                                                                                                \frac{1}{\mu_1}(Dv - d_y) \\
                                                                                                                                \frac{1}{\mu_2}(Pu - v) \\
                                                                                                                              \end{array}
                                                                                                                            \right).
\end{equation*}

\begin{algorithm}
\caption{ADAL (Anisotropic TV Denoising)}
\begin{algorithmic}[1]\label{alg:adal_tvdn}
\STATE Given $\xsupk{u}{0},\xsupk{v}{0},\lambda,\xsupk{\gamma}{0}$.
\FOR{$k = 0,1,\cdots,K$}
    \STATE $d_x \gets \mathcal{T}(Du + \mu_1 \gamma_x, \lambda\mu_1)$ \label{line:adal_dx}
    \STATE $v \gets v(d_y,u,\gamma_y,\gamma_z)$, the solution of \eqref{eq:aniso_sys_v} \label{line:adal_v}
    \STATE $d_y \gets \mathcal{T}(Dv + \mu_1\gamma_y, \lambda\mu_1)$ \label{line:adal_dy}
    \STATE $u \gets u(d_x,v,\gamma_x,\gamma_z)$, the solution of \eqref{eq:aniso_sys_u} \label{line:adal_u}
    \STATE $\gamma \gets \gamma + \Delta$\label{line:LM-update}

\ENDFOR
\RETURN $\frac{1}{2}(u+P^Tv)$
\end{algorithmic}
\end{algorithm}

\subsection{Convergence Analysis}\label{sec:conv_aniso}
We establish the convergence of Algorithm \ref{alg:adal_tvdn} by expressing problem \eqref{eq:aniso_constr} as an instance of problem \eqref{eq:struct_lin_con} and then showing that Algorithm \ref{alg:adal_tvdn} is, in fact, an ADAL method for problem \eqref{eq:struct_lin_con}, employing two-way updates to the variables.

Define $X := \left(
               \begin{array}{c}
                 d_x \\
                 v \\
               \end{array}
             \right)
$, $Y := \left(
           \begin{array}{c}
             d_y \\
             u \\
           \end{array}
         \right)
$, $f(X) := \lambda\|d_x\|_1$, and $g(Y) := \lambda\|d_y\|_1 + \frac{1}{2}\|u-b\|^2$.  Then, we can write problem \eqref{eq:aniso_constr} in the form of problem \eqref{eq:struct_lin_con} as
\begin{eqnarray}
  \min_{X,Y} && f(X) + g(Y) \label{eq:aniso_constr_general}\\
  \nonumber s.t. && AX = BY,
\end{eqnarray}
where $A = \left(
             \begin{array}{cc}
               I & 0 \\
               0 & D \\
               0 & I \\
             \end{array}
           \right) \in \bbR^{3mn\times 2mn}
$, and $B = \left(
              \begin{array}{cc}
                0 & D \\
                I & 0 \\
                0 & P \\
              \end{array}
            \right) \in \bbR^{3mn\times 2mn}
$.

Observe that Lines \ref{line:adal_dx} and \ref{line:adal_v} of Algorithm \ref{alg:adal_tvdn} exactly solve the Lagrangian subproblem of \eqref{eq:aniso_constr_general} with respect to $X$ - the subproblem is decomposable with respect to $d_x$ and $v$.  Similarly, Lines \ref{line:adal_dy} and \ref{line:adal_u} of Algorithm \ref{alg:adal_tvdn} solve the Lagrangian subproblem with respect to $Y$ - the subproblem is decomposable with respect to $d_y$ and $u$.  The matrices $A$ and $B$ obviously have full column rank.  Hence, the convergence of Algorithm \ref{alg:adal_tvdn} follows as a result of Theorem \ref{thm:admm_conv}.

\subsection{The Isotropic Case}\label{sec:iso}
The isotropic TV denoising model differs from the anisotropic model in the definition of the TV norm.  In this case, we define $\|u\|_{TV}^{ISO} := \sum_i \sqrt{(\nabla_xu)_i^2 + (\nabla_yu)_i^2} = \sum_i\|([\nabla_xu]_i,[\nabla_yu]_i)\|$, and the optimization problem to solve is
\begin{equation}\label{eq:iso_tvdn}
    \min_u \lambda\|u\|_{TV}^{ISO} + \frac{1}{2}\|u-b\|^2.
\end{equation}
Note that $\|u\|_{TV}^{ISO}$ is the group lasso regularization on $(\nabla_xu,\nabla_yu)$, with each group consisting of $([\nabla_xu]_i,[\nabla_yu]_i)$.
We introduce the same auxiliary variables and linear constraints defining them as in the previous section, except that the constraint coupling $d_y$ and $v$ becomes
\begin{equation}\label{eq:iso_constr_v}
    d_y = P^TDv.
\end{equation}
This modification is necessary because $d_x$ and $d_y$ are now coupled by the isotropic TV norm and the order of their elements have to match, i.e. column-major with respect to the original image matrix.
The subproblem with respect to $d_x$ and $d_y$ now becomes
\begin{equation}\label{eq:iso_sub_d}
    \min_{d_x,d_y} \lambda\sum_i\|([d_x]_i,[d_y]_i)\| + \gamma_x^T(Du - d_x) + \gamma_y^T(P^TDv - d_y)
    + \frac{1}{2\mu_1}(\|Du-d_x\|^2 + \|P^TDv-d_y\|^2),
\end{equation}
which is a proximal problem associated with the group $l_{1,2}$-norm $\|d\|_{1,2}$ with $d_x \equiv \nabla_xu, d_y \equiv \nabla_yu$, where the groups are defined as above.  The solution to this subproblem is thus given by a block soft-thresholding operation \cite{friedlander_spg, qin2010efficient, combettes2009proximal},
$\mathcal{S}(\left(
\begin{array}{cc}
D & 0 \\
0 & P^TD \\
\end{array}
\right)
\left(
  \begin{array}{c}
    u \\
    v \\
  \end{array}
\right) +
\mu_1 \left(
        \begin{array}{c}
          \gamma_x \\
          \gamma_y \\
        \end{array}
      \right), \lambda\mu_1
)$, where the block soft-thresholding operator is defined blockwise as
\begin{equation*}
    \mathcal{S}(x,\lambda)_i := \max\{\|x_i\|-\lambda,0\}\frac{x_i}{\|x_i\|},
\end{equation*}
and $x_i$ is the $i$-th block of $x$, i.e. $([Du+\mu_1\gamma_x]_i,[P^TDv+\mu_1\gamma_y]_i)$ in our case.
The subproblem with respect to $u$ is the same as \eqref{eq:aniso_sys_u}, and that with respect to $v$ is
\begin{equation}\label{eq:iso_sys_v}
    \left( D^TD + \frac{\mu_1}{\mu_2}I \right)v = D^TP(d_y - \mu_1\gamma_y) + \mu_2\gamma_z + \frac{\mu_1}{\mu_2}Pu.
\end{equation}

We state the ADAL method for the isotropic TV denoising in Algorithm \ref{alg:adal_tvdn_iso}, where because of \eqref{eq:iso_constr_v}, $\Delta \equiv \left(
                                                                                                                              \begin{array}{c}
                                                                                                                                \frac{1}{\mu_1}(Du - d_x) \\
                                                                                                                                \frac{1}{\mu_1}(P^TDv - d_y) \\
                                                                                                                                \frac{1}{\mu_2}(Pu - v) \\
                                                                                                                              \end{array}
                                                                                                                            \right)$.

\begin{algorithm}
\caption{ADAL (Isotropic TV Denoising)}
\begin{algorithmic}[1]\label{alg:adal_tvdn_iso}
\STATE Given $\xsupk{u}{0},\xsupk{v}{0},\lambda,\xsupk{\gamma}{0}$.
\FOR{$k = 0,1,\cdots,k_{max}$}
    \STATE $\left(
              \begin{array}{c}
                d_x \\
                d_y \\
              \end{array}
            \right) \gets \mathcal{S}\left(
                \left(
                  \begin{array}{c}
                    Du\\
                    P^TDv \\
                  \end{array}
                \right) +
                \mu_1 \left(
                        \begin{array}{c}
                          \gamma_x \\
                          \gamma_y \\
                        \end{array}
                      \right), \lambda\mu_1\right)
    $
    \STATE $v \gets v(d_y,u,\gamma_y,\gamma_z)$, the solution of \eqref{eq:iso_sys_v}
    \STATE $u \gets u(d_x,v,\gamma_x,\gamma_z)$, the solution of \eqref{eq:aniso_sys_u}
    \STATE $\gamma \gets \gamma + \Delta$ \label{line:iso_adal_lm_update}

\ENDFOR
\RETURN $\frac{1}{2}(u+P^Tv)$
\end{algorithmic}
\end{algorithm}

Due to the non-decomposability of problem \eqref{eq:iso_sub_d} with respect to $d_x$ and $d_y$ in this case, Algorithm \ref{alg:adal_tvdn_iso} cannot be interpreted as an algorithm that employs alternating updates to two blocks of variables as in Section \ref{sec:adal_tvdn}.  Hence, the convergence analysis for the anisotropic case cannot be extended to this case in a straightforward manner.  However, our experimental results in the next section show strong indication of convergence to the optimal solution.

\subsection{Comparison with The Split Bregman Method}\label{sec:contrast_sb_tvdn}
Since the split Bregman method (Algorithm \ref{alg:splitbregman}) is equivalent to the ADAL method (Algorithm \ref{alg:adal}) \cite{tai2009augmented, esser2009applications, setzer2009split} applied to the constrained problem
\begin{eqnarray*}
  \min_{d,u} && \lambda(\|d_x\|_1 + \|d_y\|_1) + \frac{1}{2}\|u-b\|^2 \\
  \nonumber s.t. && d_x = \nabla_x u, \\
  \nonumber && d_y = \nabla_y u,
\end{eqnarray*}
it is clear that the main difference between ADAL Algorithms \ref{alg:adal_tvdn} and \ref{alg:adal_tvdn_iso} and the split Bregman method comes from the introduction of the additional variable $v = Pu$ in problem \eqref{eq:aniso_constr}.  The split Bregman subproblem with respect to $u$ (line \ref{line:sb_u} in Algorithm \ref{alg:splitbregman}) can be simplified to the linear system
\begin{equation}\label{eq:sb_lin_sys}
    \left(\mu I + (\nabla_x^T\nabla_x + \nabla_y^T\nabla_y)\right)\xsupk{u}{k+1} = \mu b + \nabla_x^T(\xsupk{d_x}{k}-\xsupk{r_x}{k}) + \nabla_y^T(\xsupk{d_y}{k} - \xsupk{r_y}{k}),
\end{equation}
whose left-hand-side matrix includes a Laplacian matrix and is strictly diagonally dominant.  Solving this linear system exactly in each iteration is relatively expensive.  Hence, one iteration of the Gauss-Seidel method is applied in \cite{goldstein2009split} to solve \eqref{eq:sb_lin_sys} approximately.  Consequently, the condition for the convergence guarantee is violated in this case.

In contrast, the subproblems with respect to $v$ and $u$ in ADAL have simpler structures and thus can be solved exactly in an efficient manner as we saw in Section \ref{sec:adal_tvdn}.  The splitting of $u$ and $v$ also leads to the establishment of the global convergence of Algorithm \ref{alg:adal_tvdn} in the anisotropic case.  We surmised that this was a better approach for the TV denoising problem; the results in the next section confirmed this.

\subsection{A Globally Convergent ADAL Method for the Isotropic TV-Model}\label{sec:iso_conv}
Let us introduce three sets of variables $\left(
          \begin{array}{c}
            d_x \\
            d_y \\
          \end{array}
        \right), v$, and $w$
as follows:
\begin{equation*}
    d_x = Du, \quad d_y = P^TDv, \quad u = w, \quad v = Pw.
\end{equation*}
The isotropic TV model then has the form of \eqref{eq:aniso_constr_general} with
\begin{equation*}
    X = \left(
          \begin{array}{c}
            d_x \\
            d_y \\
            w \\
          \end{array}
        \right), Y = \left(
                   \begin{array}{c}
                     u \\
                     v \\
                   \end{array}
                 \right), A = \left(
                                \begin{array}{ccc}
                                  I & 0 & 0 \\
                                  0 & I & 0 \\
                                  0 & 0 & I \\
                                  0 & 0 & P \\
                                \end{array}
                              \right), \textrm{and } B = \left(
                                                          \begin{array}{cc}
                                                            D & 0 \\
                                                            0 & P^TD \\
                                                            I & 0 \\
                                                            0 & I \\
                                                          \end{array}
                                                        \right).
\end{equation*}
When the augmented Lagrangian $\mathcal{L}(X,Y,\gamma) \equiv \mathcal{L}(d_x,d_y,w,u,v,\gamma_x,\gamma_y,\gamma_u,\gamma_v)$ is minimized with respect to $X$, the minimization is separable in terms of $(d_x,d_y)$ and $w$.  Similarly, the minimization with respect to $Y$ is separable in terms of $u$ and $v$.

If we use the same penalty parameter $\mu_2$ for both constraints that involve $w$, then the subproblems that one obtains for $u$, $v$, and $w$ require sovling, respectively,
\begin{eqnarray}
  \left(D^TD + \left(\mu_1 + \frac{\mu_1}{\mu_2}\right)I\right)u &=& \mu_1\left( b-\gamma_u-\frac{1}{\mu_2}w \right)+D^T(d_x-\mu_1\gamma_x), \label{eq:iso_conv_u}\\
  \left(D^TD + \frac{\mu_1}{\mu_2}I\right)v &=& \mu_1\left(\frac{1}{\mu_2}Pw - \gamma_v \right) +D^TP(d_y - \mu_1\gamma_y), \label{eq:iso_conv_v}
\end{eqnarray}
and
\begin{equation}\label{eq:iso_conv_w}
    w = \frac{1}{2}\left( u + P^Tv + \mu_2(\gamma_u + P^T\gamma_v) \right).
\end{equation}
We incorporate these procedures in Algorithm \ref{alg:adal_tvdn_iso_conv} below, where now
$\gamma \equiv \left(
                                                                                                                              \begin{array}{c}
                                                                                                                                \gamma_x \\
                                                                                                                                \gamma_y \\
                                                                                                                                \gamma_u \\
																	     \gamma_v
                                                                                                                              \end{array}
                                                                                                                            \right)$, and
$\Delta \equiv \left(
                                                                                                                              \begin{array}{c}
                                                                                                                                \frac{1}{\mu_1}(Du - d_x) \\
                                                                                                                                \frac{1}{\mu_1}(P^TDv - d_y) \\
																	     \frac{1}{\mu_2}(w-u)\\
                                                                                                                                \frac{1}{\mu_2}(Pw - v) \\
                                                                                                                              \end{array}
                                                                                                                            \right)$.
Note that both Algorithms \ref{alg:adal_tvdn_iso} and \ref{alg:adal_tvdn_iso_conv} compute $u, v, \gamma_x$, and $\gamma_y$.  In addition,  Algorithm \ref{alg:adal_tvdn_iso_conv} requires the computation of $w, \gamma_u$, and $\gamma_v$ whereas the isotropic ADAL Algorithm \ref{alg:adal_tvdn_iso} only requires computation of $\gamma_z$.  Consequently, the new convergent algorithm has slightly more work at each iteration.

\begin{algorithm}
\caption{ADAL (Isotropic TV Denoising) - Convergent}
\begin{algorithmic}[1]\label{alg:adal_tvdn_iso_conv}
\STATE Given $\xsupk{u}{0},\xsupk{v}{0},\lambda,\xsupk{\gamma}{0}$.
\FOR{$k = 0,1,\cdots,k_{max}$}
    \STATE $\left(
              \begin{array}{c}
                d_x \\
                d_y \\
              \end{array}
            \right) \gets \mathcal{S}\left(
                \left(
                  \begin{array}{c}
                    Du \\
                    P^TDv \\
                  \end{array}
                \right) +
                \mu_1 \left(
                        \begin{array}{c}
                          \gamma_x \\
                          \gamma_y \\
                        \end{array}
                      \right), \lambda\mu_1\right)
    $
    \STATE $w \gets w(u,v,\gamma_v,\gamma_u$, by \eqref{eq:iso_conv_w}
    \STATE $v \gets v(d_y,w,\gamma_y,\gamma_v)$, the solution of \eqref{eq:iso_conv_v}
    \STATE $u \gets u(d_x,w,\gamma_x,\gamma_u)$, the solution of \eqref{eq:iso_conv_u}
    \STATE $\gamma \gets \gamma + \Delta$

\ENDFOR
\RETURN $\frac{1}{3}(u+P^Tv+w)$
\end{algorithmic}
\end{algorithm}

\subsection{Practical Implementation}
It is often beneficial to associate a step-size $\theta$ to the Lagrange multiplier updates, i.e. $\gamma \gets \gamma + \theta\Delta$ at line \ref{line:LM-update} in Algorithms \ref{alg:adal_tvdn} and \ref{alg:adal_tvdn_iso_conv} and line \ref{line:iso_adal_lm_update} in Algorithm \ref{alg:adal_tvdn_iso}.  The convergence of this ADAL variant with $\theta \in (0,\frac{\sqrt{5}+1}{2})$ has been established in \cite{glowinski1989augmented},\cite{he2000alternating}, and \cite{wen2010alternating} under various contexts.  In our implementation, we set $\theta = 1.618$.

In practice, we can often set $\mu_1 = \mu_2 = \mu$.  In this case, we can save some scalar-vector multiplications by maintaining $\tilde{\gamma}=\mu\gamma$ instead of the $\gamma$ variables themselves.  Then, for example, \eqref{eq:iso_conv_u} can be simplified to
\begin{equation*}
    \left(D^TD + \left(\mu + 1\right)I\right)u = \mu b - \tilde{\gamma_u} - w + D^T(d_x - \tilde{\gamma_x}).
\end{equation*}
All the steps in the proposed algorithms remain the same with this substitution, except for the computation of $\Delta$, which no longer involves division by $\mu$.

In addition, if $\mu$ is kept constant throughout the algorithms, we can factorize the constant left-hand-sides of the tri-diagonal linear systems and cache the factors for subsequent iterations.  The factors are lower-triangular and are stored in the form of two vectors representing the diagonal and the sub-diagonal of these factors.  Then, we can solve for $u$ and $v$ quickly through forward and backward substitutions, which require $5mn$ flops each.  In our implementation, we used the LAPACK routines \textbf{dpttrf} (for factorization) and \textbf{dpttrs} (for forward/backward substitutions).

We also developed a simple updating scheme for $\mu$, which decreases $\mu$ by a constant factor $\kappa$ after every $J$ iterations, starting from $\bar{\mu}$, bounded below by $\underline{\mu}$, i.e.
\begin{equation*}
    \xsupk{\mu}{k} = \max(\underline{\mu}, \frac{\bar{\mu}}{\kappa^{\frac{k}{J}}}).
\end{equation*}
Such an updating scheme allows different values of $\mu$ to be applied at different stages of convergence.  From our computational experience in the next Section, this often led to improved convergence speed.  We denote this variant by ``ADAL-$\mu$" and ``ADAL-conv-$\mu$" corresponding to ADAL (Algorithms \ref{alg:adal_tvdn} and \ref{alg:adal_tvdn_iso}) and Algorithm \ref{alg:adal_tvdn_iso_conv} respectively.

\section{Experiments}\label{sec:exp}
We implemented our ADAL algorithms (Algorithms \ref{alg:adal_tvdn},\ref{alg:adal_tvdn_iso},and \ref{alg:adal_tvdn_iso_conv}) in C++ with BLAS and LAPACK routines.  SplitBregman is in C with a Matlab interface. \footnote{Code downloaded from \emph{http://www.stanford.edu/~tagoldst/code.html}.}  We ran all the algorithms on a laptop with an Intel Core i5 Duo processor and 6G memory.

\subsection{Test Images}
\begin{figure}
\begin{center}
    \hspace*{-0.0in}\includegraphics[width=1.0\textwidth]{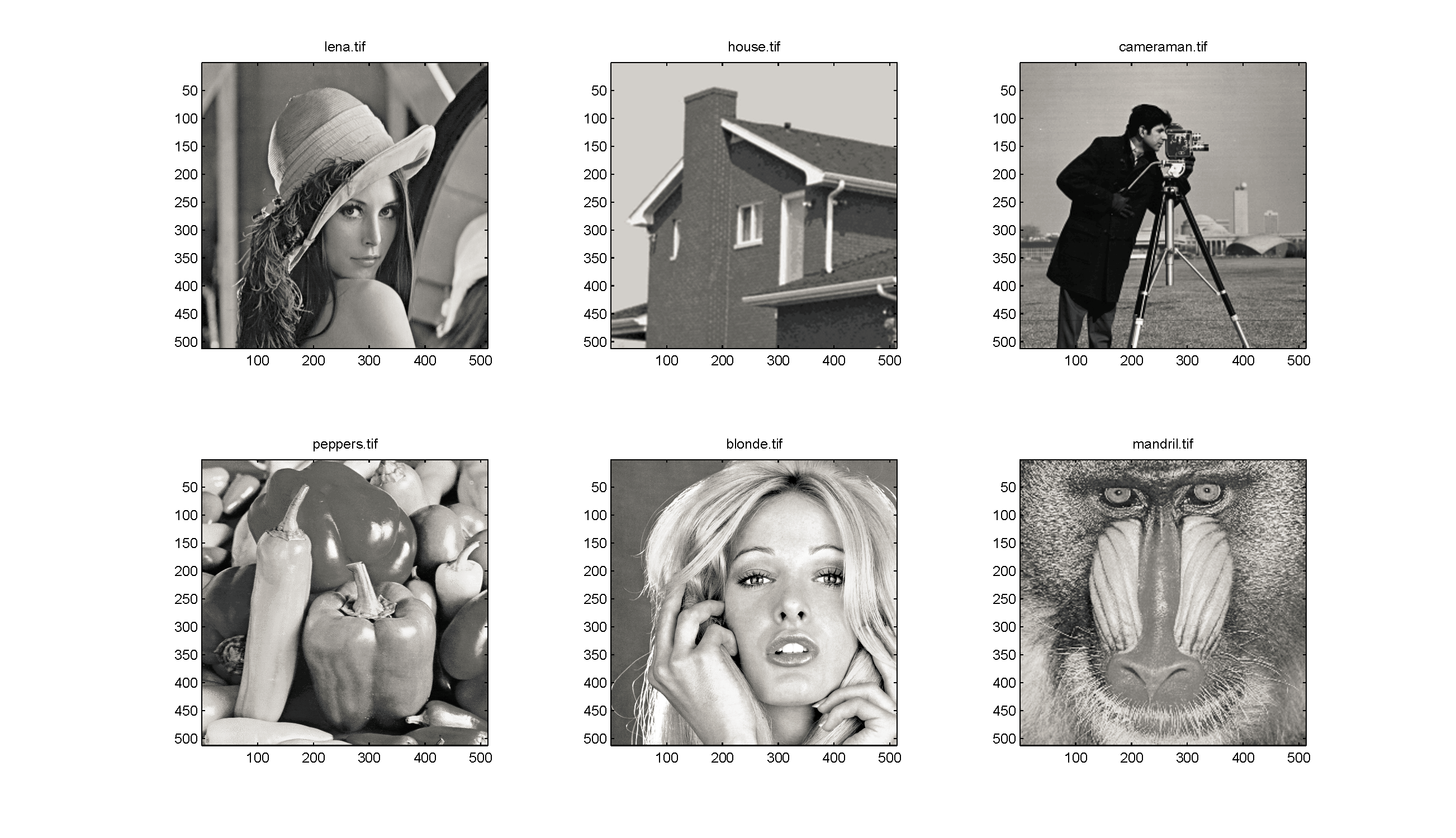}
    \caption{The set of standard test images.}
    \label{fig:test_images}
\end{center}
\end{figure}

We compared our ADAL algorithm with the split Bregman method on a set of six standard test images: \textbf{lena}, \textbf{house}, \textbf{cameraman}, \textbf{peppers}, \textbf{blonde}, and \textbf{mandril} (see Figure \ref{fig:test_images}).  They present a range of challenges to image denoising algorithms, such as the reproduction of fine detail and textures, sharp transitions and edges, and uniform regions.  Each image is a $512\times 512$ array of grey-scale pixels and is denoted by $u_0$ in vectorized form.

\subsection{Set-up}
We constructed noisy images by adding Gaussian noise to the original images, i.e. $b = u_0 + \epsilon$, where $\epsilon \sim \mathcal{N}(0,\sigma^2)$ and $b$ is the vectorized noisy image.  We set $\sigma = 30$, which introduced a considerable amount of noise.  The quality of the denoised image in the $k$-th iteration, $\usupk$ is measured by the normalized error with respect to a high quality reference solution $u^*$, i.e. $\eta^{(k)} = \frac{\|\usupk-u^*\|}{\|u^*\|}$, as well as the peak-signal-to-noise ratio (PSNR). The PSNR of an image $u$ with respect to the noiseless image $u_0$, in the case where the maximum pixel magnitude is 255, is defined as
\begin{equation*}
    PSNR = 20\log_{10}\left( \frac{255\sqrt{nm}}{\|u-u_0\|} \right).
\end{equation*}
PSNR is monotone decreasing with the $\|u-u_0\|$, i.e. a higher PSNR indicates better reconstruction quality.

In practice, the algorithms can be stopped once an acceptable level of optimality has been reached.  For ADAL, we used the maximum of the relative primal and dual residuals \cite{boyd2010distributed}, denoted by $\epsilon$ to approximately measure the optimality of the solution.
For each image, we computed a reference solution $u^*$ and the corresponding PSNR $p^*$ by running Algorithm \ref{alg:adal_tvdn} for the anisotropic TV model and algorithm \ref{alg:adal_tvdn_iso_conv} for the isotropic TV model until the measure $\epsilon$ fell below $10^{-12}$.
We then recorded the number of iterations $K$ required by each of the algorithms to reach a normalized error $\xsupk{\eta}{K}$ less than $10^{-5}$.  We also recorded the number of iterations required to reach a PSNR $p$ whose relative gap to $p^*$, gap $\equiv \frac{|p^*-p|}{p^*}$, was less than $10^{-3}$.

\begin{figure}
\begin{center}
    \hspace*{-0.7in}\includegraphics[width=0.4\textwidth]{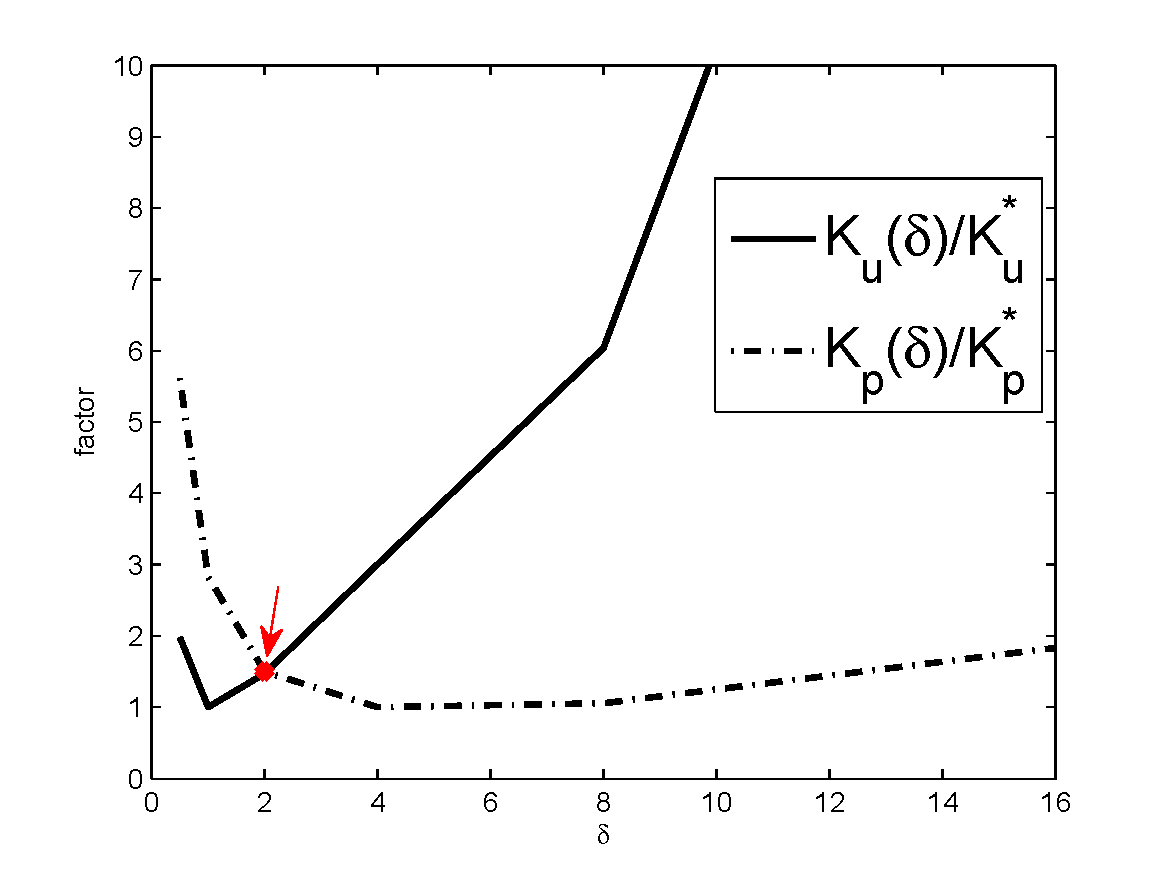}
    \caption{Plots of ratios $\frac{K_p(\delta)}{K_p^*}$ and $\frac{K_u(\delta)}{K_u^*}$ as functions of $\delta$ for isotropic ADAL on the image \textbf{house}.  The optimal point on the frontier is labeled in red.}
    \label{fig:mu_selection_adal}
\end{center}
\end{figure}

We set all initial values ($\xsupk{u}{0}, \xsupk{v}{0},\xsupk{w}{0},\xsupk{\gamma}{0}$) to zeros.  We tuned all the algorithms under investigation for convergence speed with respect to both $u^*$ and $p^*$ to reach the tolerance levels defined above, and the same set of parameters were used for all six images.  For ADAL, for both the anisotropic and isotropic TV models, we set $\mu_1 = \mu_2$ so that there was only one parameter to tune.  We tried values of $\mu$ from the set $\{0.1\delta: \delta=\frac{1}{2},1,2,4,8,16\}$ and found the value $\mu = 0.2$ to work best; the results reported below used this value.  Figure \ref{fig:mu_selection_adal} illustrates this selection using the image \textbf{house} as an example, by plotting the ratios $\frac{K_p(\delta)}{K_p^*}$ and $\frac{K_u(\delta)}{K_u^*}$ as functions of $\delta$ for isotropic ADAL, where $K_u(\delta)$ is the number of iterations needed to reduce the normalized error $\xsupk{\eta}{K}$ below the tolerance $10^{-5}$ with $\mu=\delta 0.1$, and $K_u^* = \min_\delta\{K_u(\delta)\}$.  $K_p(\delta)$ and $K_p^*$ are defined similarly with respect to $p$.   For ADAL-$\mu$ and ADAL-conv-$\mu$, we set $\kappa=1.5, J = 50, \bar{\mu}=0.5$, and $\underline{\mu}=0.05$.  For SplitBregman, we set $\mu = \frac{4}{\lambda}$, which was the value of $\mu$ from the set $\{\frac{\delta}{\lambda} : \delta=\frac{1}{2},1,2,4,8,16\}$, chosen in the same manner as that for ADAL.  In general, the parameters produced fairly consistent performance on the set of different images.

%

\subsection{Convergence Comparisons}
\begin{figure}
    \hspace*{-0.7in}\includegraphics[width=1.2\textwidth]{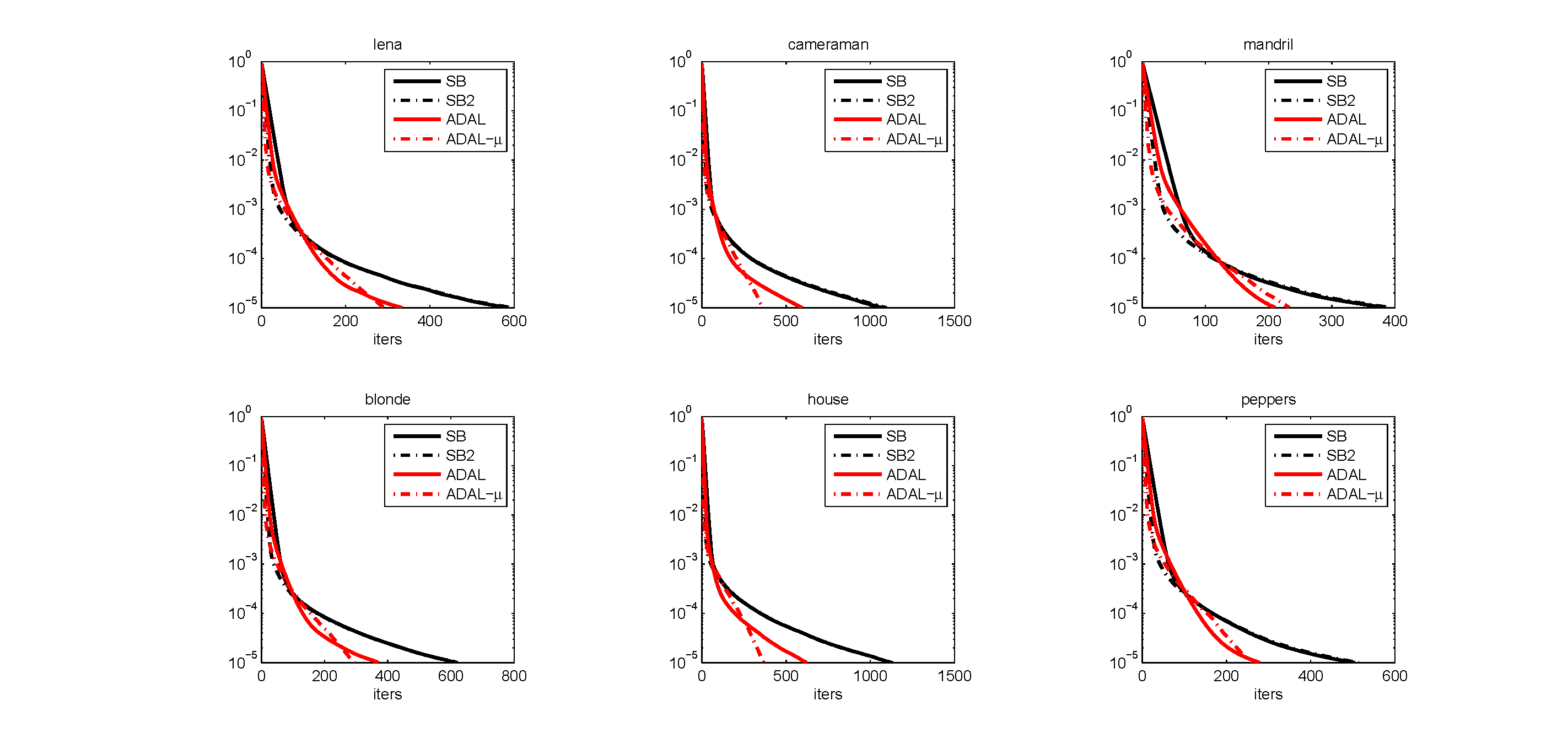}
    \caption{Convergence plots of normalized errors w.r.t. the reference solution against iterations for the anisotropic TV model.}
    \label{fig:conv_errs}
\end{figure}

\begin{figure}
    \hspace*{-0.7in}\includegraphics[width=1.2\textwidth]{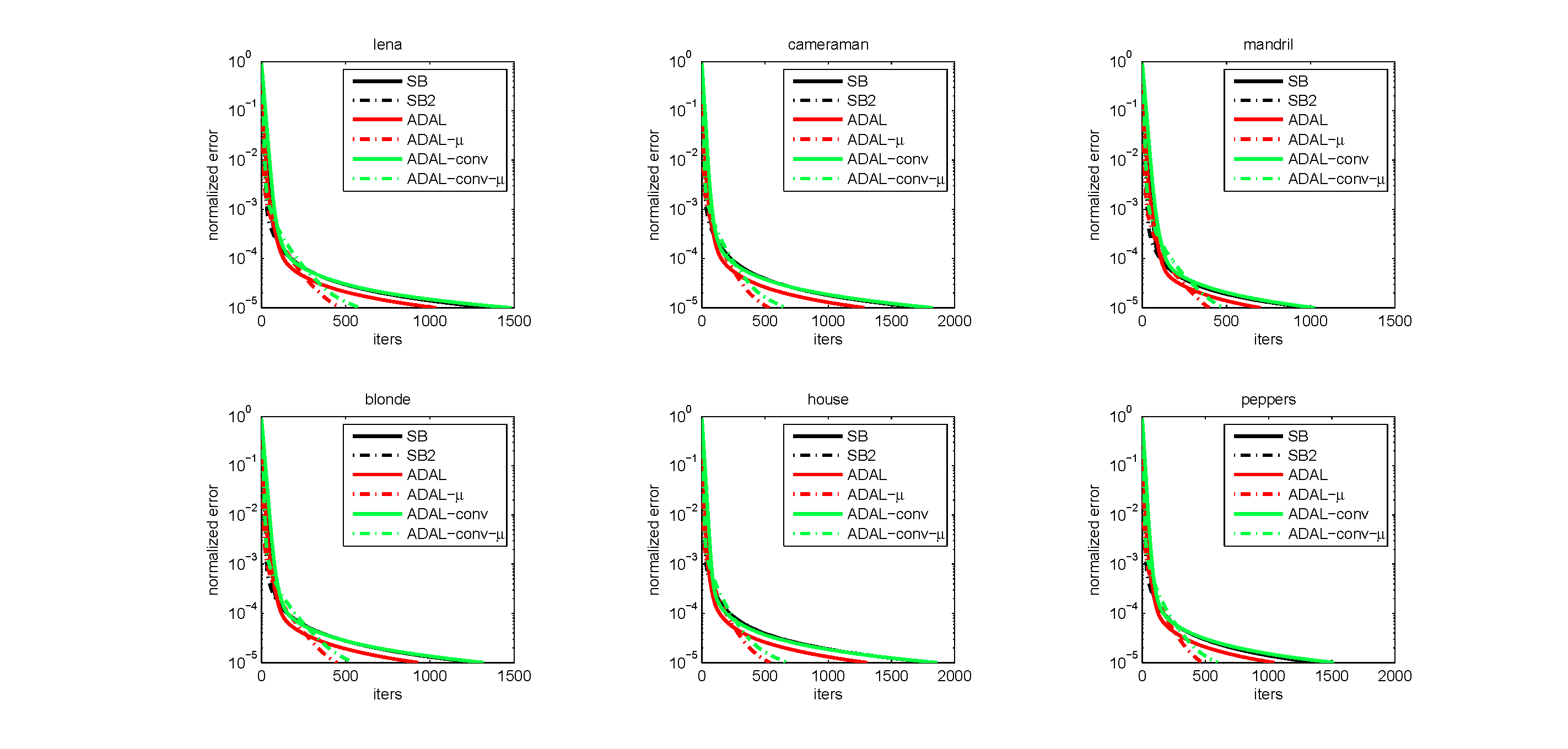}
    \caption{Convergence plots of normalized errors w.r.t. the reference solution against iterations for the isotropic TV model.}
    \label{fig:conv_errs_iso}
\end{figure}

We now present experimental results for both the anisotropic and the isotropic TV models. ``ADAL-conv'' denotes the convergent ADAL Algorithm \ref{alg:adal_tvdn_iso_conv} for the isotropic model.  We also tested a version of SplitBregman, denoted by ``SplitBregman2", where two cycles of Gauss-Seidel were performed to solve the linear system \eqref{eq:sb_lin_sys}.  In Table \ref{tab:iters}, we report the number of iterations and the corresponding CPU time required by the three algorithms to reach the stopping criterion discussed above.  Figures \ref{fig:conv_errs} and \ref{fig:conv_errs_iso} plot the relative gaps with respect to the reference solution as a function of the iteration number for all the algorithms.

In general, ADAL required fewer iterations than SplitBregman to reach the prescribed relative gap with respect to the reference solution.  The difference was particularly significant at high accuracy levels, about $45\%$ for the anisotropic model and $25\%$ for the isotropic model.  We believe that ADAL benefits from the fact that it is able to solve its subproblems exactly and efficiently, while the approximation made in the solution to the linear system in the iterations of SplitBregman slows down the convergence as its iterates approach the exact solution.

\begin{figure}
    \hspace*{-0.7in}\includegraphics[width=1.2\textwidth]{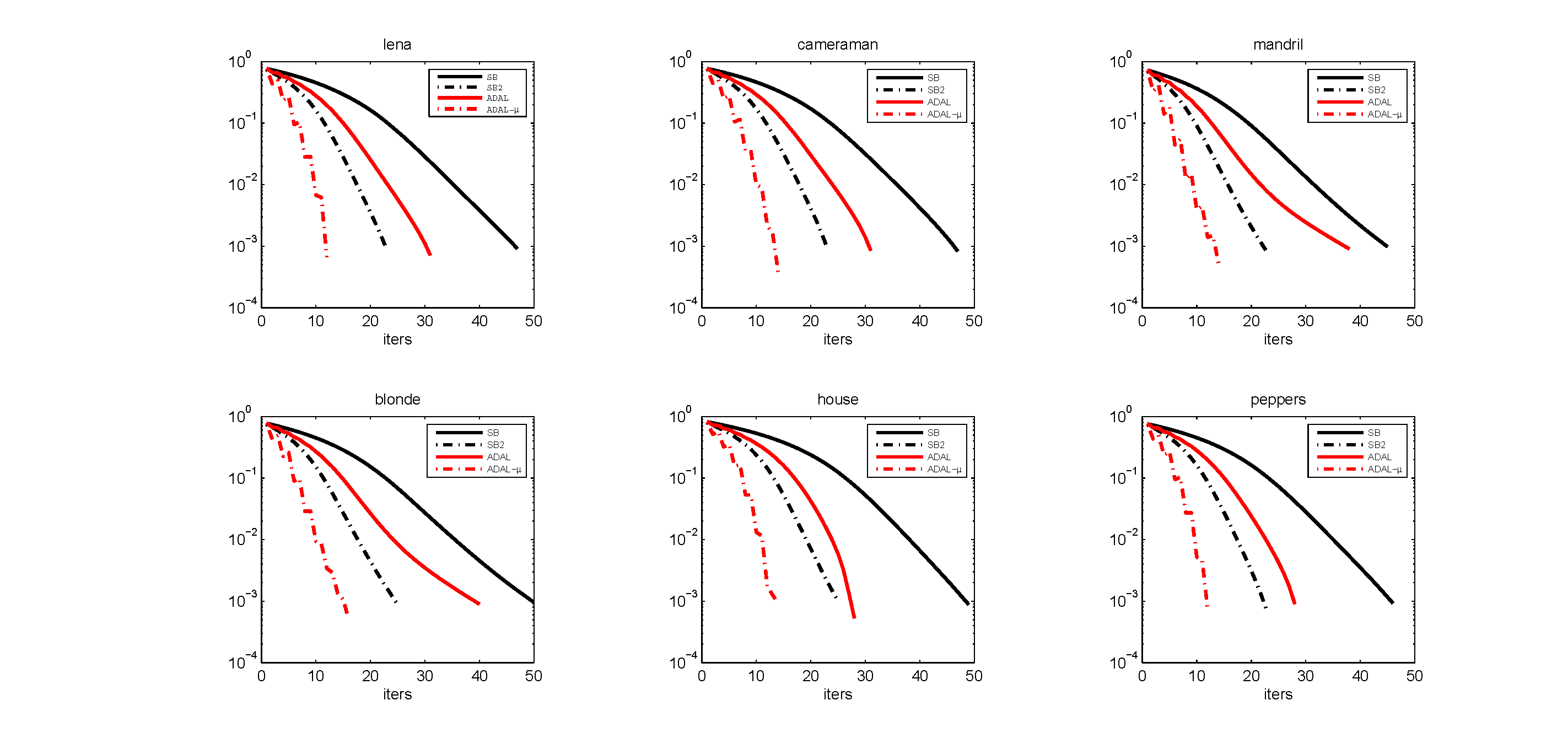}
    \caption{Convergence plots of normalized errors w.r.t. the reference PSNR against iterations for the anisotropic TV model.}
    \label{fig:psnr_errs}
\end{figure}

\begin{figure}
    \hspace*{-0.7in}\includegraphics[width=1.2\textwidth]{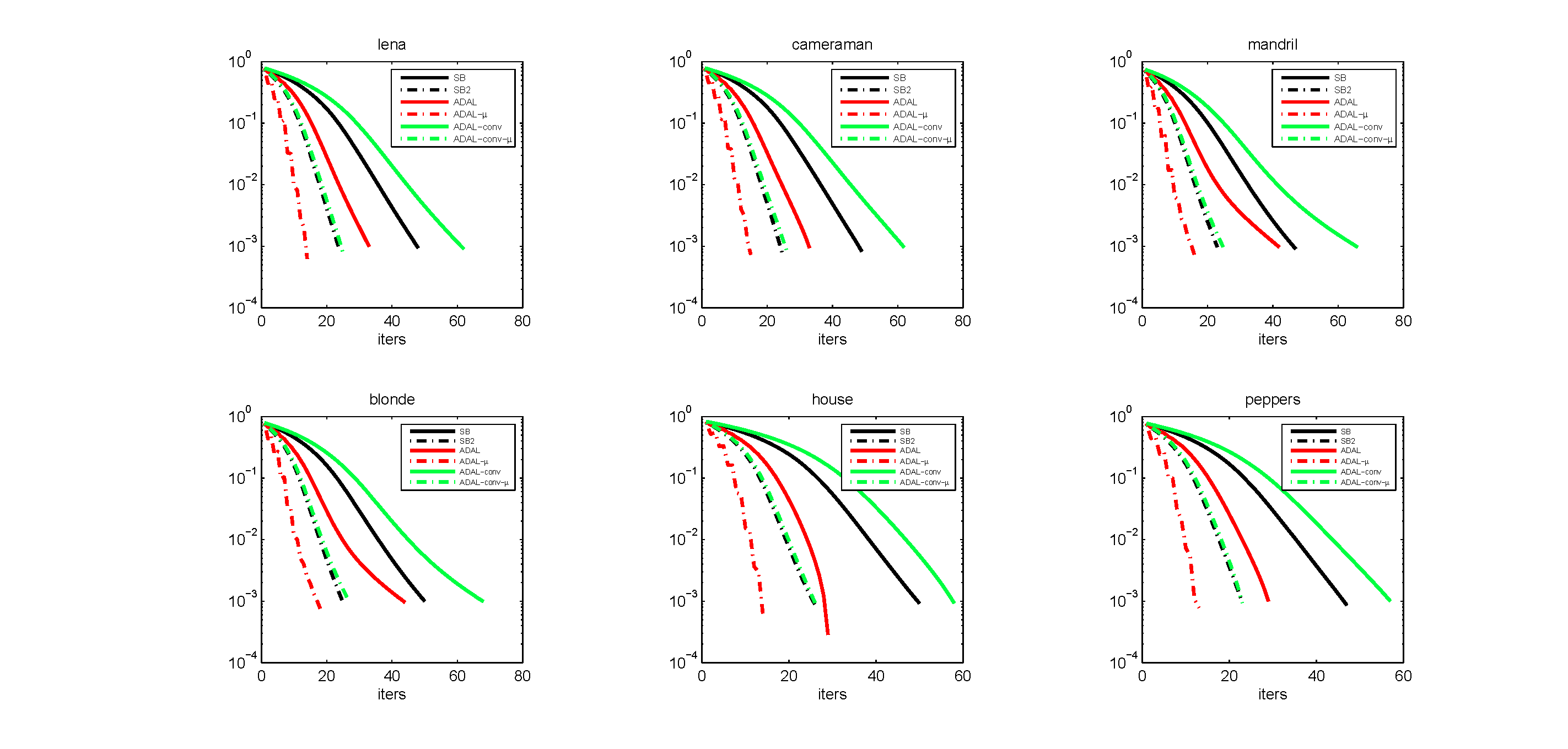}
    \caption{Convergence plots of normalized errors w.r.t. the reference PSNR against iterations for the isotropic TV model.}
    \label{fig:psnr_errs_iso}
\end{figure}

Figures \ref{fig:psnr_errs} and \ref{fig:psnr_errs_iso} plot the relative gaps with respect to the reference PSNR as a function of the iteration number.  In terms of denoising speed, ADAL also requried fewer iterations than SplitBregman to reach the denoising quality of the reference solution for both TV models.

SplitBregman2 generally required about the same number of iterations as SplitBregman to reach within the gap tolerance with respect to $u^*$, while requiring half the number of iterations as SplitBregman to reach the prescribed gap with respect to $p^*$.  It appeared that the additional cycle of Gauss-Seidel in each iteration helped improve the convergence speed initially but not in the later stages.

We  note that ADAL was faster than ADAL-conv for the isotropic model.  This appears to be due to the need for ADAL-conv to make the norm of the difference between an additional set of variables and what they are defined to be equal to small enough so as to obtain high quality images.

\begin{table}
\begin{center}
\hspace*{-0.45in}
\begin{tabular}{|c|c|c|c|c|c|c|c|c|c|c|}
\hline\hline
\multirow{2}{*}{Model} & \multirow{2}{*}{Algs} & \multicolumn{3}{c|}{\textbf{lena}} & \multicolumn{3}{c|}{\textbf{cameraman}}	& \multicolumn{3}{c|}{\textbf{mandril}}  \\
\hhline{~~---------}
& & iters & P-iters & CPU &  iters & P-iters & CPU &  iters & P-iters & CPU \\
\hline
\multirow{2}{*}{Anisotropic}
& ADAL & 334 & 31 & 11.1 	& 595  &31 & 18.9	& 210 &  38 & 7.0	 \\
& ADAL-$\mu$  & 291 & 11 & 9.2 	& 360 & 14& 11.5	& 232& 14& 7.4\\
& SplitBregman &584 & 47&  12.5	& 1070 & 47& 23.6	& 385& 45 & 8.6	  	\\
& SplitBregman2 &595 & 23&  15.7	& 1093 & 23& 29.0	& 394& 23 & 10.8	  	\\
\hline
\multirow{2}{*}{Isotropic}
& ADAL & 1036 &33 & 37.2 	&1286 &33 &46.0 	&701 &42 & 25.2	   	\\
& ADAL-$\mu$  & 472 & 14& 17.1 	& 531& 15& 19.3	& 396& 16& 14.5\\
& ADAL-conv 	  & 1482 &62 & 54.8 	& 1826& 62& 67.0	 & 1023& 66& 38.2   \\
& ADAL-conv-$\mu$  & 587 &25 &22.1  	&666 &26 &25.2 	& 482& 25& 18.3\\
& SplitBregman & 1372& 48& 32.4 	& 1767& 49& 41.6	&983 &47 & 23.2	   	\\
& SplitBregman2 & 1376& 24& 41.7 	& 1779& 25& 53.5	&996 &23 & 29.4	   	\\
\hline\hline
\end{tabular}

\vspace*{0.2in}
\hspace*{-0.45in}
\begin{tabular}{|c|c|c|c|c|c|c|c|c|c|c|}
\hline\hline
\multirow{2}{*}{Model} & \multirow{2}{*}{Algs} & \multicolumn{3}{c|}{\textbf{blonde}}	& \multicolumn{3}{c|}{\textbf{house}}	& \multicolumn{3}{c|}{\textbf{peppers}}  \\
\hhline{~~---------}
& & iters & P-iters & CPU &  iters & P-iters & CPU &  iters & P-iters & CPU \\
\hline
\multirow{2}{*}{Anisotropic}
& ADAL & 370 & 40 & 12.0 	&621 &28 & 20.7	& 279& 28& 9.4	 \\
& ADAL-$\mu$  & 293 &15 & 9.3 	& 364& 14& 11.6	& 262& 12& 8.4\\
& SplitBregman & 617& 50& 13.6 	& 1126&49 &25.6 	&504 &46 &11.2 	  	\\
& SplitBregman2 & 619& 25& 17.2 	& 1128&25 &33.1 	&513 &23 &14.0 	  	\\
\hline
\multirow{2}{*}{Isotropic}
& ADAL &929 & 44& 34.1 	&1308 &29 &49.4 	&1043 &29 & 43.4   	\\
& ADAL-$\mu$  & 447& 18& 16.4 	&537 & 14& 19.5	& 476& 13& 18.9 \\
& ADAL-conv & 1320& 68& 51.5 	& 1860& 58& 70.6	& 1518& 57& 68.1	     \\
& ADAL-conv-$\mu$  & 547 &27 & 20.6 	& 677& 26& 25.6	&595 & 23& 26.5\\
& SplitBregman & 1292& 50& 30.4 	& 1848& 50& 42.9	&1368 & 47& 33.7	   	\\
& SplitBregman2 & 1298& 25& 37.8 	& 1856& 25& 55.6	&1375 & 23& 42.5	   	\\
\hline\hline
\end{tabular}
\end{center}
\caption{Computational statistics.  ``iters" denotes the number of iterations to reach within a gap of 1e-5 w.r.t. the reference solution.  ``P-iters" denotes number of iterations to reach within a gap of 1e-3 w.r.t. the reference PSNR.  CPU time is in seconds and corresponds to the ``iters" column.}
\label{tab:iters}
\end{table}

\begin{figure}
    \hspace*{-0.55in}\vspace*{-0.0in}\includegraphics[width=0.6\textwidth]{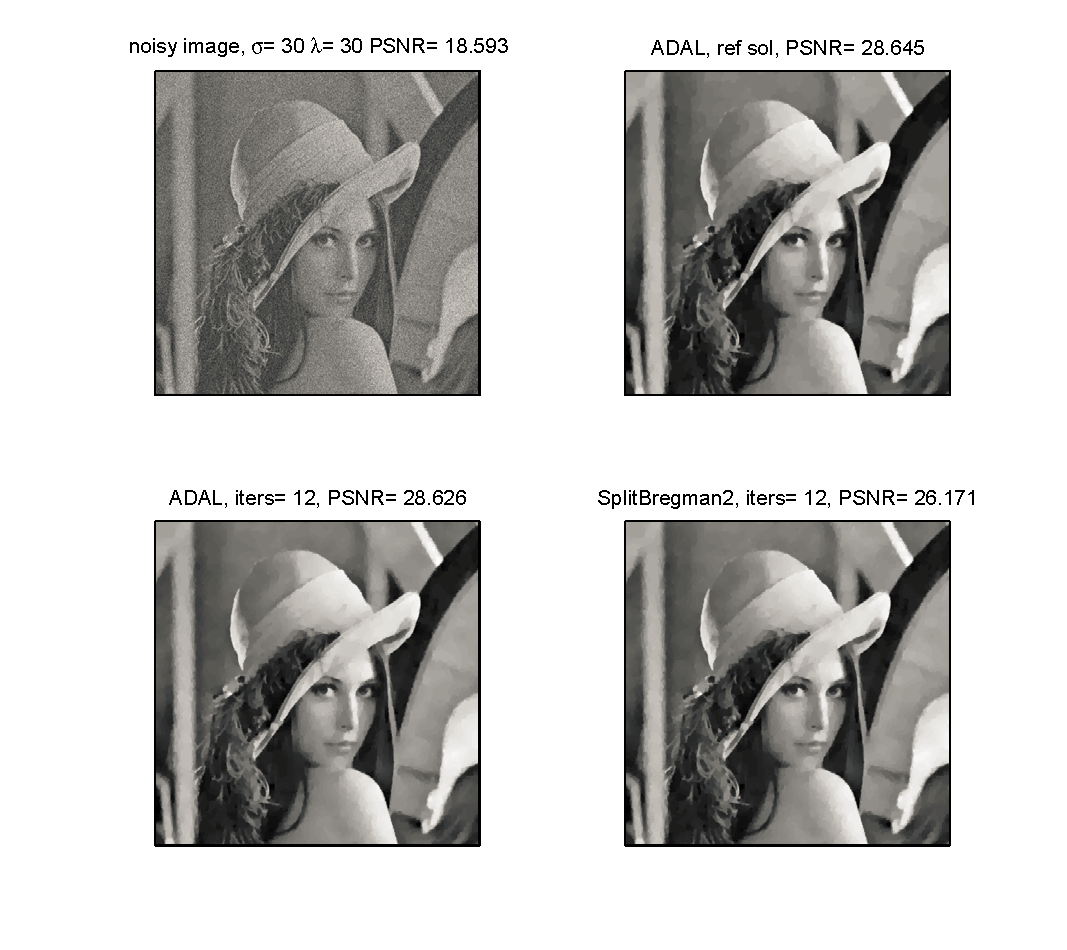}
    \hspace*{-0.4in}\vspace{-0.0in}\includegraphics[width=0.6\textwidth]{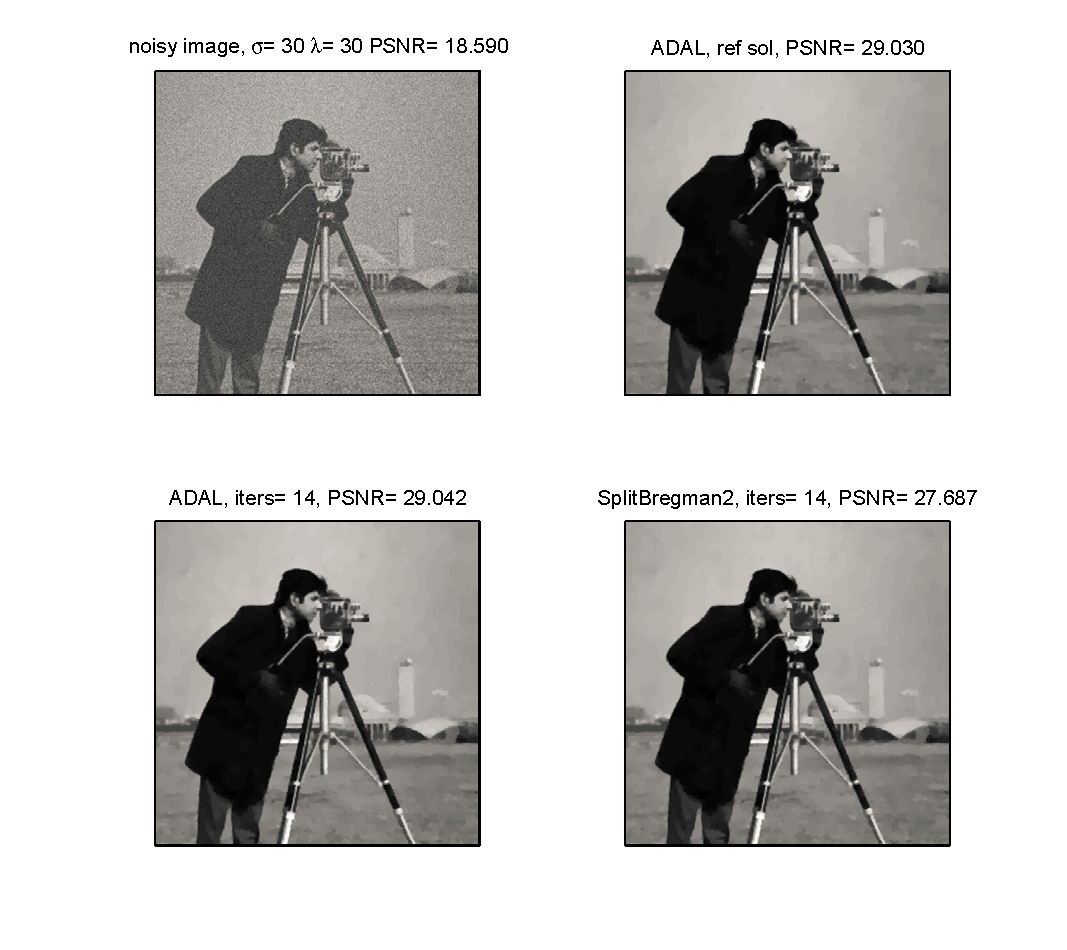}
    \hspace*{-0.55in}\vspace*{-0.0in}\includegraphics[width=0.6\textwidth]{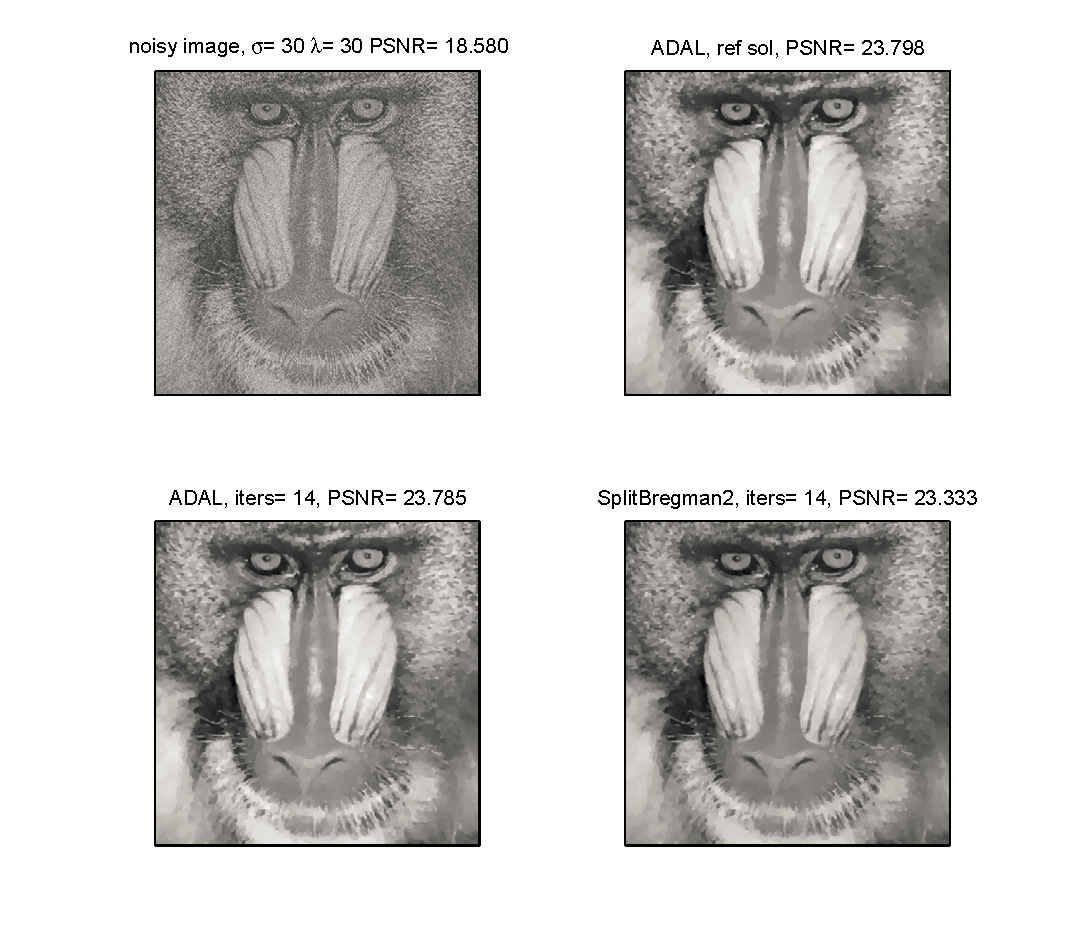}
    \hspace*{-0.4in}\vspace*{-0.0in}\includegraphics[width=0.6\textwidth]{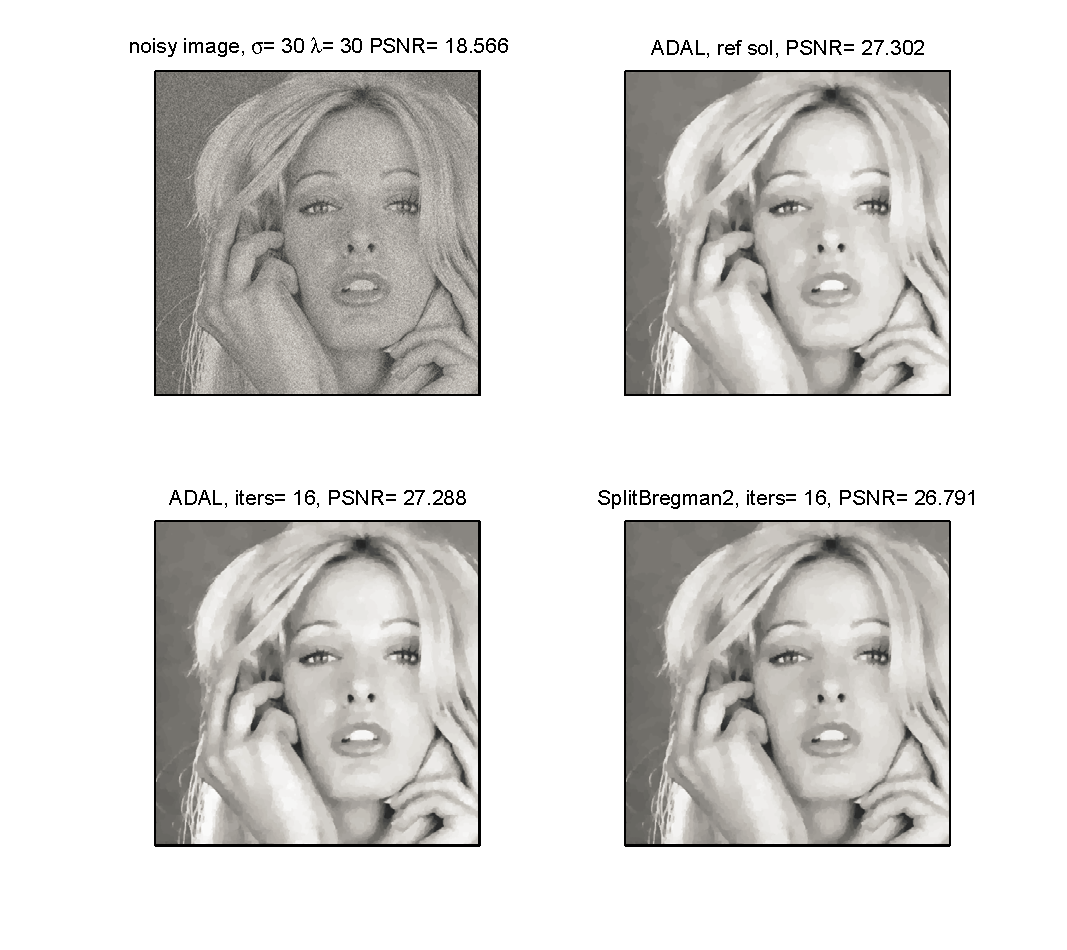}
    \caption{Comparison of reconstruction quality for \textbf{lena}, \textbf{cameraman}, \textbf{mandril}, and \textbf{blonde} with the anisotropic TV model.  Top left: noisy image. Top right: reference solution obtained by ADAL. Bottom left: ADAL solution obtained after the corresponding number of iterations indicated in Table \ref{tab:iters}. Bottom right: SplitBregman solution obtained after the same number of iterations.}
    \label{fig:quality_1}
\end{figure}

\begin{figure}
    \hspace*{-0.55in}\vspace*{-0.0in}\includegraphics[width=0.6\textwidth]{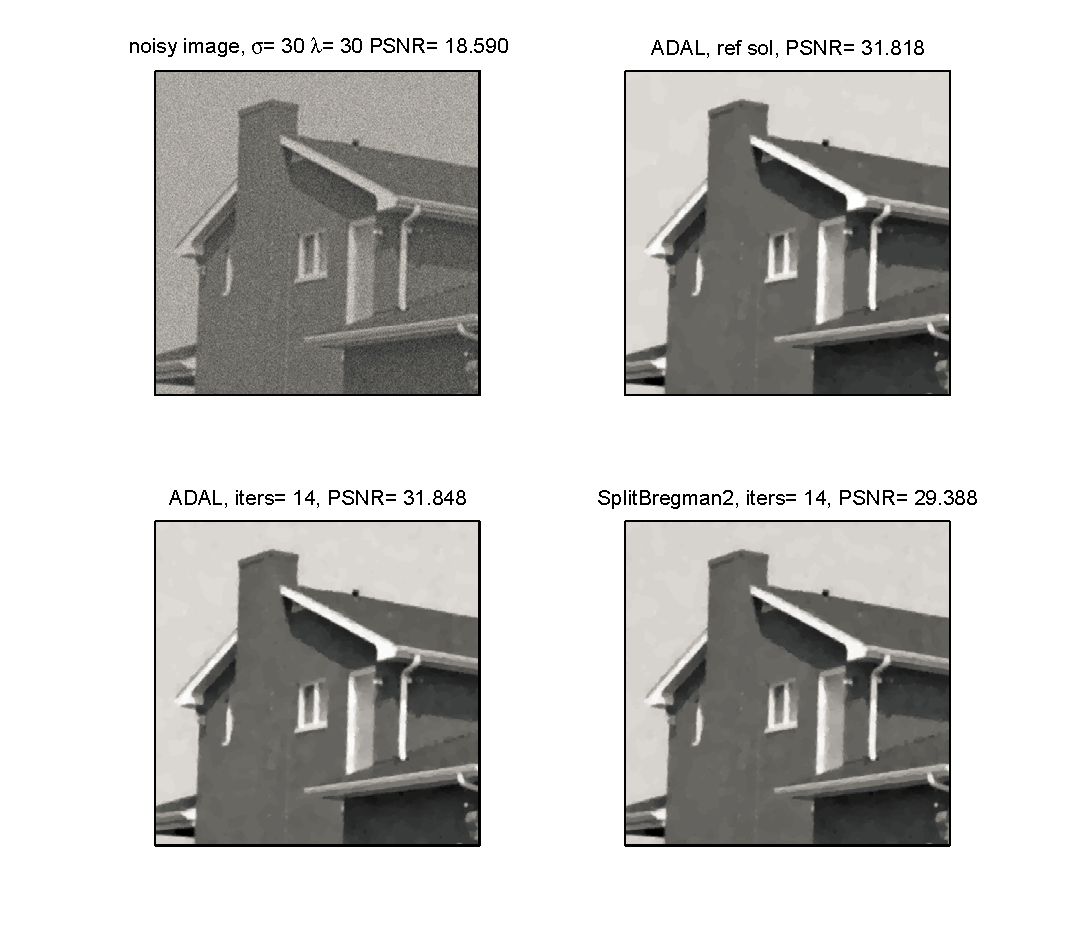}
    \hspace*{-0.4in}\vspace{-0.0in}\includegraphics[width=0.6\textwidth]{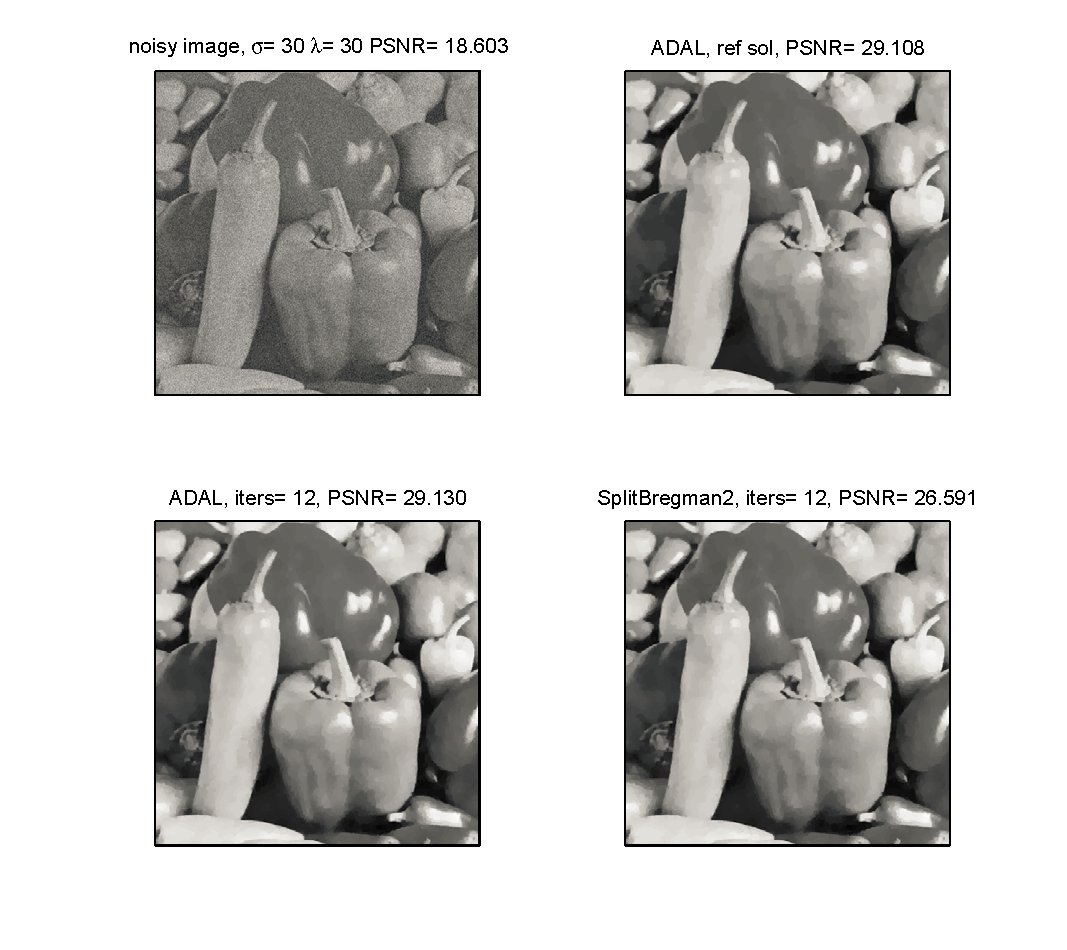}
    \caption{Comparison of reconstruction quality for \textbf{house}, and \textbf{peppers} with the anisotropic TV model.  Top left: noisy image. Top right: reference solution obtained by ADAL. Bottom left: ADAL-$\mu$ solution obtained after the corresponding number of iterations indicated in Table \ref{tab:iters}. Bottom right: SplitBregman solution obtained after the same number of iterations.}
    \label{fig:quality_2}
\end{figure}

\begin{figure}
    \hspace*{-0.55in}\vspace*{-0.0in}\includegraphics[width=0.6\textwidth]{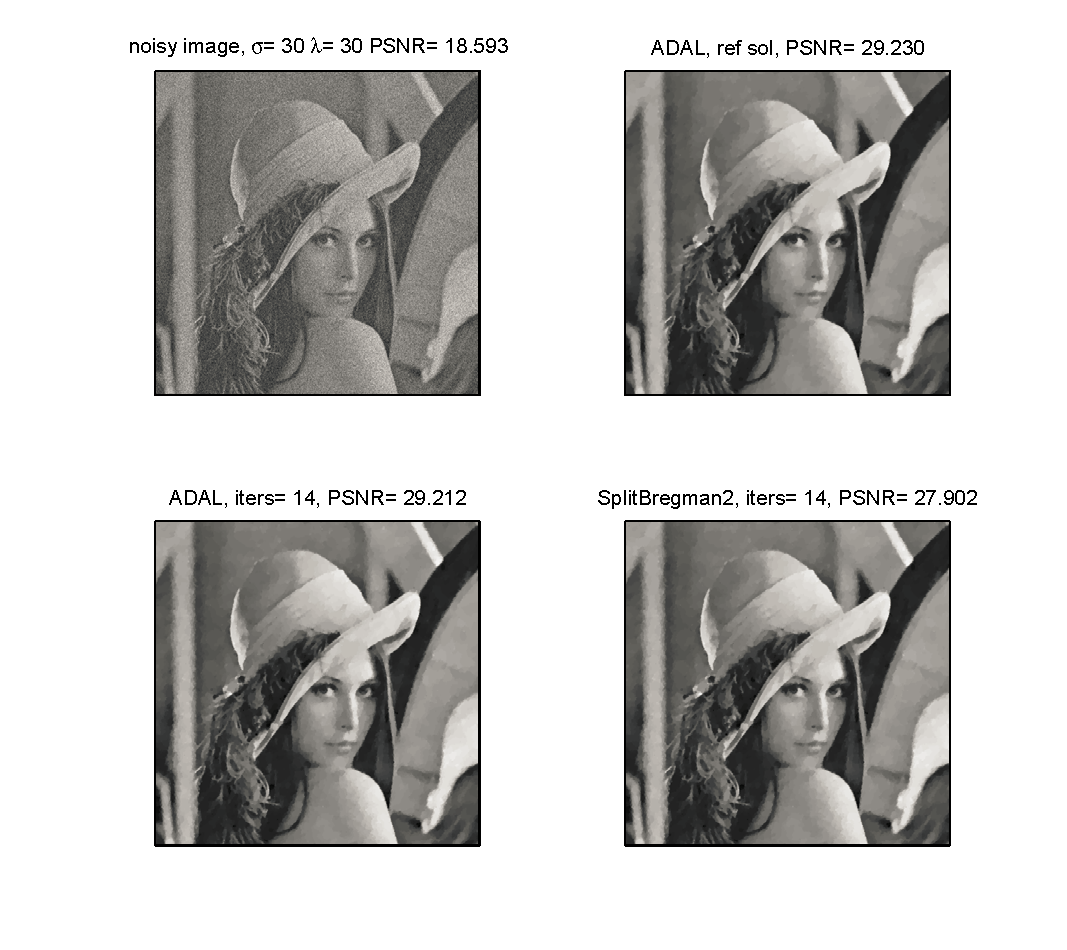}
    \hspace*{-0.4in}\vspace{-0.0in}\includegraphics[width=0.6\textwidth]{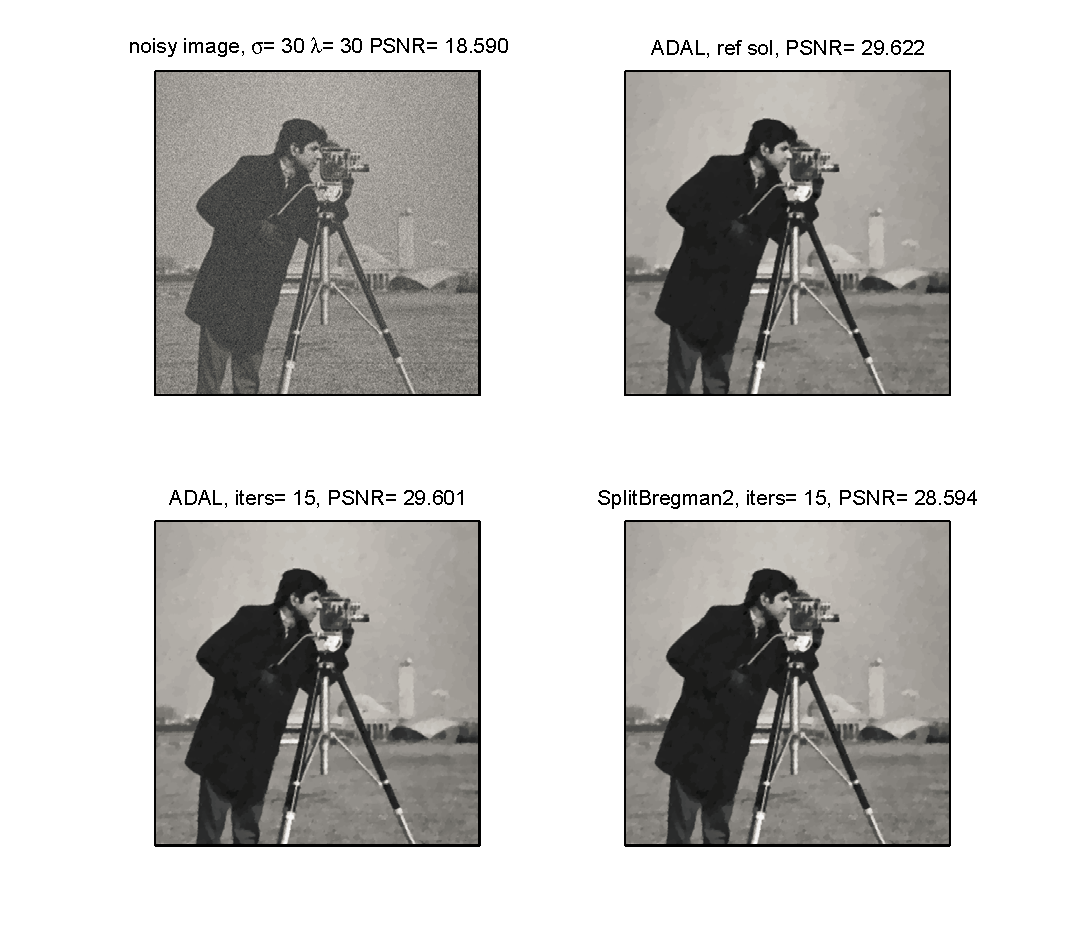}
    \hspace*{-0.55in}\vspace*{-0.0in}\includegraphics[width=0.6\textwidth]{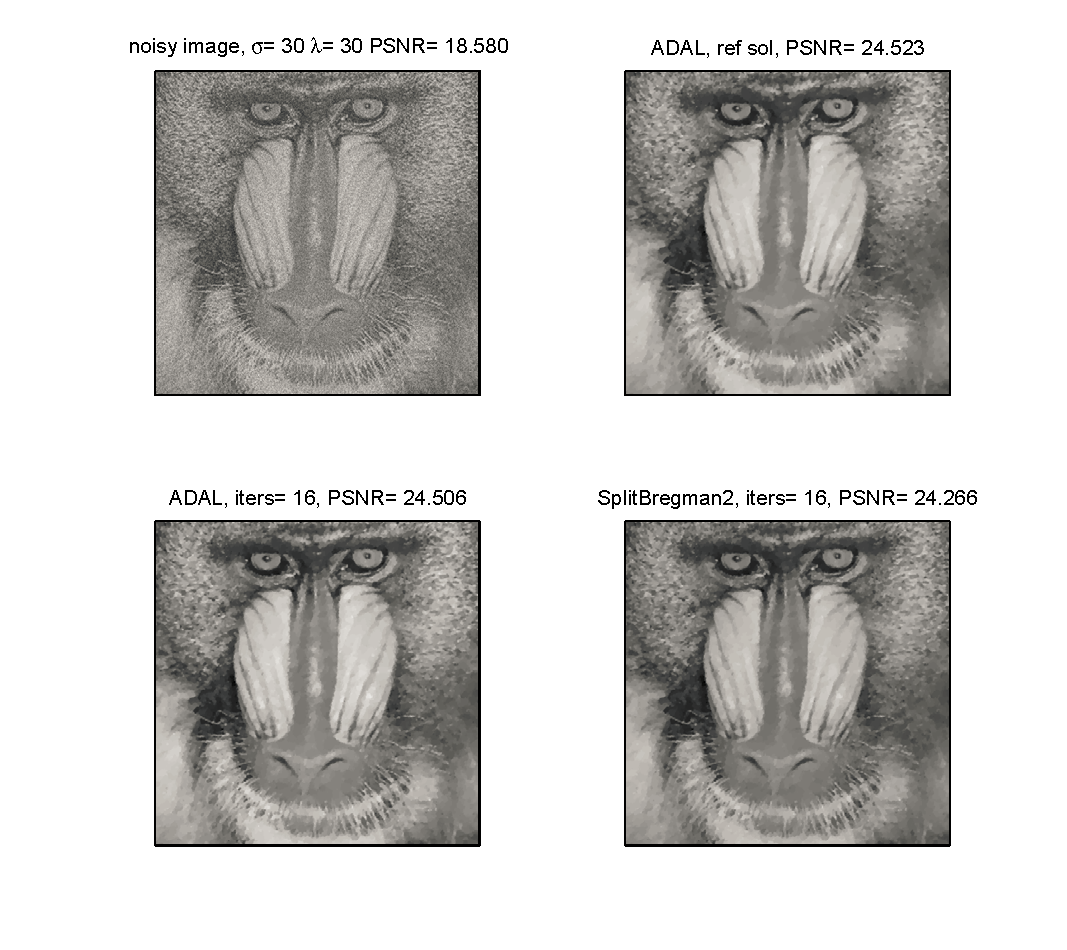}
    \hspace*{-0.4in}\vspace*{-0.0in}\includegraphics[width=0.6\textwidth]{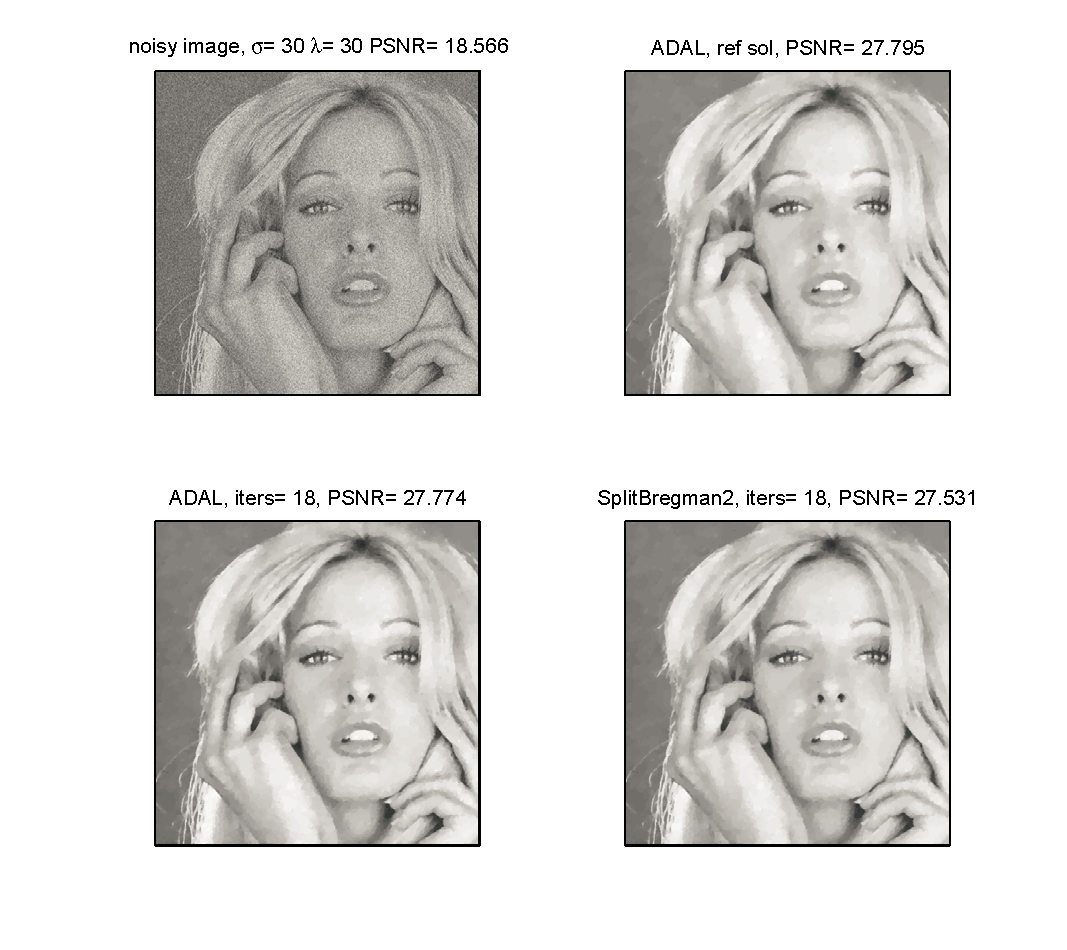}
    \caption{Comparison of reconstruction quality for \textbf{lena}, \textbf{cameraman}, \textbf{mandril}, and \textbf{blonde} with the isotropic TV model.  Top left: noisy image. Top right: reference solution obtained by ADAL. Bottom left: ADAL-$\mu$ solution obtained after the corresponding number of iterations indicated in Table \ref{tab:iters}. Bottom right: SplitBregman solution obtained after the same number of iterations.}
    \label{fig:quality_1_iso}
\end{figure}

\begin{figure}
    \hspace*{-0.55in}\vspace*{-0.0in}\includegraphics[width=0.6\textwidth]{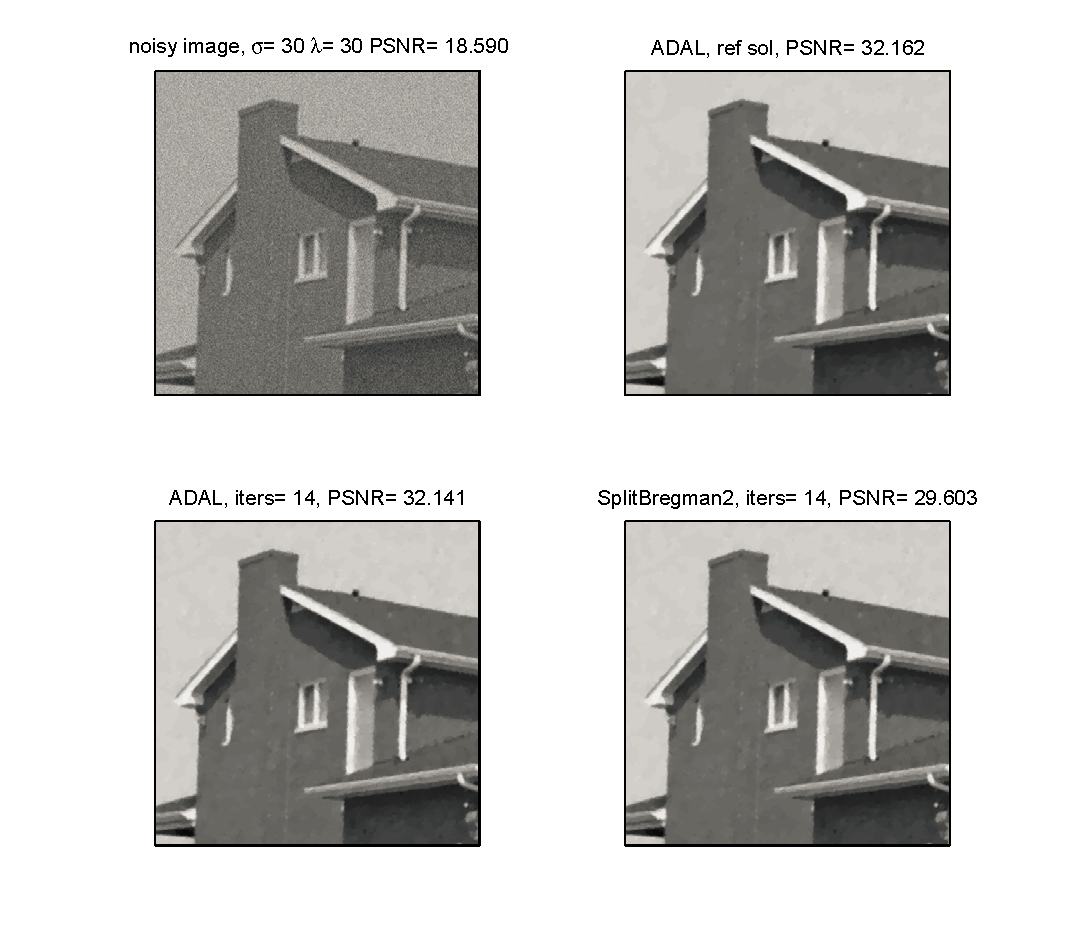}
    \hspace*{-0.4in}\vspace{-0.0in}\includegraphics[width=0.6\textwidth]{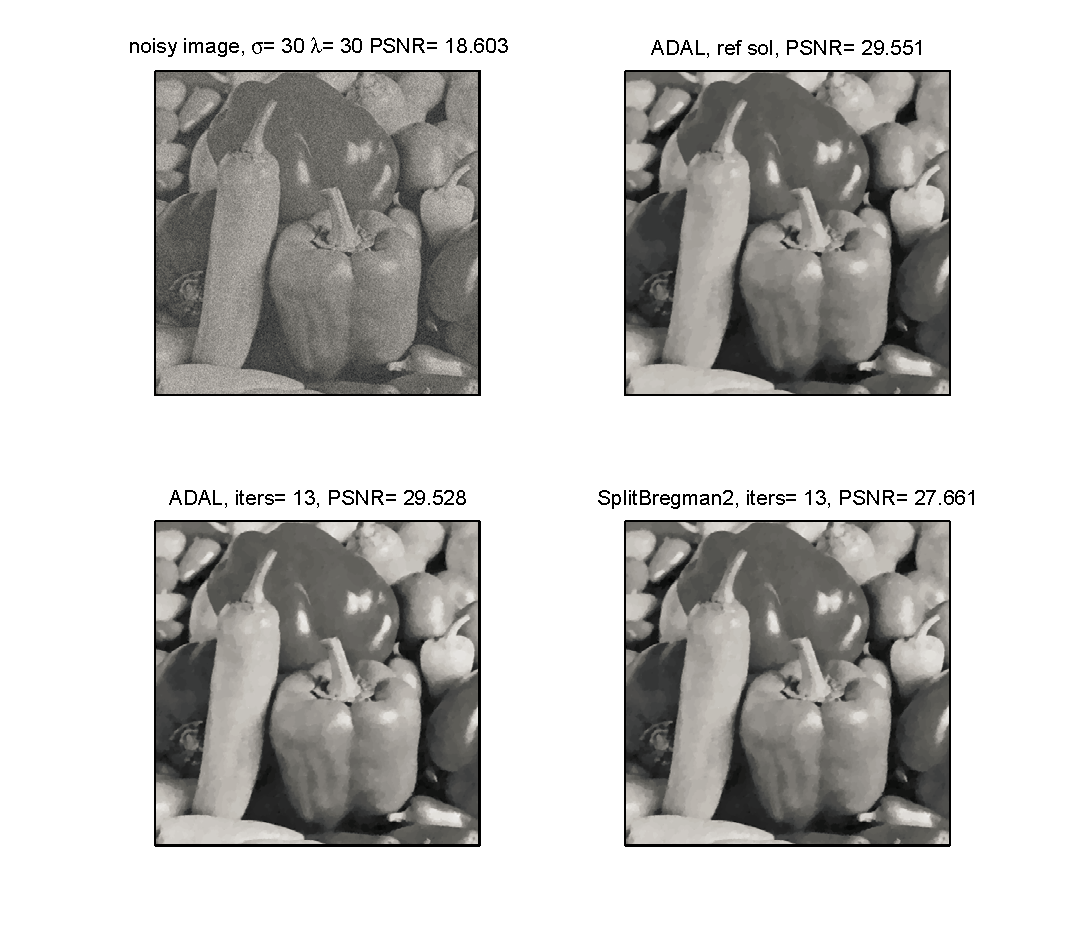}
    \caption{Comparison of reconstruction quality for \textbf{house}, and \textbf{peppers} with the isotropic TV model.  Top left: noisy image. Top right: reference solution obtained by ADAL. Bottom left: ADAL solution obtained after the corresponding number of iterations indicated in Table \ref{tab:iters}. Bottom right: SplitBregman solution obtained after the same number of iterations.}
    \label{fig:quality_2_iso}
\end{figure}

We observe that the updating scheme for $\mu$ improved the speed of convergence, especially for the isotropic TV model.  ADAL-$\mu$ and ADAL-conv-$\mu$ reduced the number of iterations required by ADAL and ADAL-conv, to achieve the same normalized errors on the six test images by an amount between $40\%$ and $60\%$, and $60\%$ and $67\%$, respectivelly.  We show in Figures \ref{fig:quality_1}, \ref{fig:quality_1_iso}, \ref{fig:quality_2}, and \ref{fig:quality_2_iso} the solutions obtained by ADAL-$\mu$ and SplitBregman after the number of iterations specified in the ADAL-$\mu$ row in Table \ref{tab:iters}.

In terms of CPU time, all of the algorithms required more time to reach the prescribed relative gap with respect to the reference solution for the isotropic model than that for the anisotropic model because of more iterations required.  For the anisotropic model, ADAL took less time than SplitBregman, and ADAL-$\mu$ took the least, with the least number of iterations.  The apparent more per-iteration work for ADAL is due to the computation for the additional tridiagonal linear system and the Lagrange multiplier, as compared to one sweep of Gauss-Seidel employed by SplitBregman.  SplitBregman2 required longer time than SplitBregman because of the additional sweep of Gaus-Seidel per iteration with about the same number of iterations.  The approximate total number of flops per iteration\footnote{Depending on the implementation, the number of flops required may vary slightly.} (including solving the linear systems, updating the Lagrange multipliers, and performing the shrinkage operations) is 33 $mn$, 46 $mn$, and 44 $mn$ for SplitBregman, SplitBregman2, and ADAL respectively.

For the isotropic model, ADAL-$\mu$ and ADAL-conv-$\mu$ required the least amount of time to reach the same level of accuracy, with about half the number of iterations required by SplitBregman.   On the other hand, ADAL took more time than SplitBregman in this case, though still requiring less time than SplitBregman2.  The driving factor here was the greater number of iterations required than in the anisotropic case.  ADAL-conv required slightly more per-iteration work than ADAL due to the additional computation for $w$ and the Lagrange multiplier for the additional constraint introduced.  The approximate total number of flops per iteration is 34 $mn$, 47 $mn$, 48 $mn$, and 57 $mn$ for SplitBregman, SplitBregman2, ADAL, and ADAL-conv respectively.


\section{Conclusion}
We have proposed new ADAL algorithms for solving TV denoising problems in image processing.  The key feature of our algorithms is their use of multiple variable splittings which results in their ability to solve the ADAL subproblems exactly and efficiently.  Our first ADAL algorithm has a global convergence guarantee for the case of anisotropic TV model, and the experimental results show that its iterates converge significantly faster than those of SplitBregman.  Even though the convergence guarantee of this ADAL variant cannot be extended easily to the isotropic TV model, empirical results show that with a simple updating scheme for $\mu$, it still compares favorably to SplitBregman in convergence speed.  We also proposed another ADAL variant for the isotropic TV model that has a global convergence guarantee, with a slightly higher per-iteration computational cost.  Because of the additional variable splitting required to obtain the convergence guarantee, the method also takes more iterations than the simpler isotropic ADAL variant.

\section{Acknowledgement}
We would like to thank the anonymous referees for their careful reading and valuable comments, which have helped improve the paper significantly. This work was supported in part by DMS 10-16571, ONR Grant N00014-08-1-1118 and DOE Grant DE-FG02-08ER25856.  Research of Shiqian Ma was supported in part by a Direct Grant of the Chinese University of Hong Kong (Project ID: 4055016) and the Hong Kong Research Grants Council Early Career Scheme (Project ID: CUHK 439513).

\bibliographystyle{abbrv}
\bibliography{tony_bib}

\begin{thebibliography}{10}

\bibitem{afonso2010fast}
M.~Afonso, J.~Bioucas-Dias, and M.~Figueiredo.
\newblock Fast image recovery using variable splitting and constrained
  optimization.
\newblock {\em Image Processing, IEEE Transactions on}, 19(9):2345--2356, 2010.

\bibitem{afonso2009augmented}
M.~Afonso, J.~Bioucas-Dias, and M.~Figueiredo.
\newblock {An augmented Lagrangian approach to the constrained optimization
  formulation of imaging inverse problems}.
\newblock {\em IEEE Transactions on Image Processing}, (20):681--695, 2011.

\bibitem{almeida2012deconvolving}
M.~S. Almeida and M.~A. Figueiredo.
\newblock Deconvolving images with unknown boundaries using the alternating
  direction method of multipliers.
\newblock {\em IEEE Transactions on Image Processing, to appear (arXiv
  1210.02687v2)}.

\bibitem{bertsekas1999nonlinear}
D.~Bertsekas.
\newblock {\em {Nonlinear Programming}}.
\newblock Athena Scientific Belmont, MA, 1999.

\bibitem{boyd2010distributed}
S.~Boyd, N.~Parikh, E.~Chu, B.~Peleato, and J.~Eckstein.
\newblock Distributed optimization and statistical learning via the alternating
  direction method of multipliers.
\newblock {\em Foundations and Trends in Machine Learning}, 3(1):1--123, 2010.

\bibitem{chambolle2004algorithm}
A.~Chambolle.
\newblock An algorithm for total variation minimization and applications.
\newblock {\em Journal of Mathematical Imaging and Vision}, 20(1):89--97, 2004.

\bibitem{combettes2009proximal}
P.~Combettes and J.~Pesquet.
\newblock Proximal splitting methods in signal processing.
\newblock {\em Fixed-Point Algorithms for Inverse Problems in Science and
  Engineering}, pages 185--212, 2011.

\bibitem{deng2011group}
W.~Deng, W.~Yin, and Z.~Y.
\newblock Group sparse optimization by alternating direction method.
\newblock Technical report, TR 11-06, Rice University, 2011.

\bibitem{eckstein1992douglas}
J.~Eckstein and D.~Bertsekas.
\newblock {On the Douglas-Rachford splitting method and the proximal point
  algorithm for maximal monotone operators}.
\newblock {\em Mathematical Programming}, 55(1):293--318, 1992.

\bibitem{esser2009applications}
E.~Esser.
\newblock Applications of lagrangian-based alternating direction methods and
  connections to split bregman.
\newblock {\em CAM report}, 9:31, 2009.

\bibitem{figueiredo2010restoration}
M.~A. Figueiredo and J.~M. Bioucas-Dias.
\newblock Restoration of poissonian images using alternating direction
  optimization.
\newblock {\em Image Processing, IEEE Transactions on}, 19(12):3133--3145,
  2010.

\bibitem{gabay1976dual}
D.~Gabay and B.~Mercier.
\newblock {A dual algorithm for the solution of nonlinear variational problems
  via finite element approximation}.
\newblock {\em Computers \& Mathematics with Applications}, 2(1):17--40, 1976.

\bibitem{glowinski1989augmented}
R.~Glowinski and P.~Le~Tallec.
\newblock {\em Augmented Lagrangian and operator-splitting methods in nonlinear
  mechanics}, volume~9.
\newblock SIAM, 1989.

\bibitem{glowinski1975adal}
R.~Glowinski and A.~Marroco.
\newblock {Sur l'approximation, par elements finis d'ordre un, et la
  resolution, par penalisation-dualite d'une classe de problemes de dirichlet
  non lineares}.
\newblock {\em Rev. Francaise d'Automat. Inf. Recherche Operationelle},
  (9):41--76, 1975.

\bibitem{goldfarb2009parametric}
D.~Goldfarb and W.~Yin.
\newblock Parametric maximum flow algorithms for fast total variation
  minimization.
\newblock {\em SIAM Journal on Scientific Computing}, 31(5):3712--3743, 2009.

\bibitem{goldstein2009split}
T.~Goldstein and S.~Osher.
\newblock The split bregman method for l1-regularized problems.
\newblock {\em SIAM Journal on Imaging Sciences}, 2:323, 2009.

\bibitem{golub1996matrix}
G.~Golub and C.~Van~Loan.
\newblock {\em {Matrix Computations}}.
\newblock Johns Hopkins Univ Pr, 1996.

\bibitem{he2002new}
B.~He, L.~Liao, D.~Han, and H.~Yang.
\newblock A new inexact alternating directions method for monotone variational
  inequalities.
\newblock {\em Mathematical Programming}, 92(1):103--118, 2002.

\bibitem{he2000alternating}
B.~He, H.~Yang, and S.~Wang.
\newblock Alternating direction method with self-adaptive penalty parameters
  for monotone variational inequalities.
\newblock {\em Journal of Optimization Theory and applications},
  106(2):337--356, 2000.

\bibitem{hestenes1969multiplier}
M.~Hestenes.
\newblock Multiplier and gradient methods.
\newblock {\em Journal of Optimization Theory and Applications}, 4(5):303--320,
  1969.

\bibitem{lin2010augmented}
Z.~Lin, M.~Chen, L.~Wu, and Y.~Ma.
\newblock {The augmented lagrange multiplier method for exact recovery of
  corrupted low-rank matrices}.
\newblock {\em Arxiv Preprint arXiv:1009.5055}, 2010.

\bibitem{ma2009fixed}
S.~Ma, D.~Goldfarb, and L.~Chen.
\newblock Fixed point and bregman iterative methods for matrix rank
  minimization.
\newblock {\em Mathematical Programming}, pages 1--33, 2009.

\bibitem{nocedal}
J.~Nocedal and S.~Wright.
\newblock {\em {Numerical Optimization}}.
\newblock Springer Verlag, 1999.

\bibitem{osher2006iterative}
S.~Osher, M.~Burger, D.~Goldfarb, J.~Xu, and W.~Yin.
\newblock An iterative regularization method for total variation-based image
  restoration.
\newblock {\em Multiscale Modeling and Simulation}, 4(2):460--489, 2006.

\bibitem{powell1972nonlinear}
M.~Powell.
\newblock A method for nonlinear constraints in minimization problems.
\newblock In R.~Fletcher, editor, {\em Optimization}. Academic Press, New York,
  New York, 1972.

\bibitem{qin2011structured}
Z.~Qin and D.~Goldfarb.
\newblock Structured sparsity via alternating direction methods.
\newblock {\em Journal of Machine Learning Research}, 13:1373--1406, 2012.

\bibitem{qin2010efficient}
Z.~Qin, K.~Scheinberg, and D.~Goldfarb.
\newblock Efficient block-coordinate descent algorithms for the group lasso.
\newblock {\em Mathematical Programming Computation}, 5:143--169, 2013.

\bibitem{ramani2012splitting}
S.~Ramani and J.~A. Fessler.
\newblock A splitting-based iterative algorithm for accelerated statistical
  x-ray ct reconstruction.
\newblock {\em IEEE Transactions on Medical Imaging}, 31(3):677--688, 2012.

\bibitem{rockafellar1973multiplier}
R.~Rockafellar.
\newblock The multiplier method of hestenes and powell applied to convex
  programming.
\newblock {\em Journal of Optimization Theory and Applications},
  12(6):555--562, 1973.

\bibitem{rudin1992nonlinear}
L.~Rudin, S.~Osher, and E.~Fatemi.
\newblock Nonlinear total variation based noise removal algorithms.
\newblock {\em Physica D: Nonlinear Phenomena}, 60(1-4):259--268, 1992.

\bibitem{setzer2009split}
S.~Setzer.
\newblock Split bregman algorithm, douglas-rachford splitting and frame
  shrinkage.
\newblock {\em Scale space and variational methods in computer vision}, pages
  464--476, 2009.

\bibitem{shen2011augmented}
Y.~Shen, Z.~Wen, and Y.~Zhang.
\newblock Augmented lagrangian alternating direction method for matrix
  separation based on low-rank factorization.
\newblock {\em TR11-02, Rice University}, 2011.

\bibitem{steidl2010removing}
G.~Steidl and T.~Teuber.
\newblock Removing multiplicative noise by douglas-rachford splitting methods.
\newblock {\em Journal of Mathematical Imaging and Vision}, 36(2):168--184,
  2010.

\bibitem{strong2003edge}
D.~Strong and T.~Chan.
\newblock Edge-preserving and scale-dependent properties of total variation
  regularization.
\newblock {\em Inverse problems}, 19:S165, 2003.

\bibitem{tai2009augmented}
X.~Tai and C.~Wu.
\newblock Augmented lagrangian method, dual methods and split bregman iteration
  for rof model.
\newblock {\em Scale Space and Variational Methods in Computer Vision}, pages
  502--513, 2009.

\bibitem{tao2009alternating}
M.~Tao and J.~Yang.
\newblock Alternating direction algorithms for total variation deconvolution in
  image reconstruction.
\newblock {\em Optimization Online}, 2009.

\bibitem{friedlander_spg}
E.~van~den Berg, M.~Schmidt, M.~Friedlander, and K.~Murphy.
\newblock {Group sparsity via linear-time projection}.
\newblock Technical report, TR-2008-09, Department of Computer Science,
  University of British Columbia, 2008.

\bibitem{wang2008new}
Y.~Wang, J.~Yang, W.~Yin, and Y.~Zhang.
\newblock A new alternating minimization algorithm for total variation image
  reconstruction.
\newblock {\em SIAM Journal on Imaging Sciences}, 1(3):248--272, 2008.

\bibitem{wen2010alternating}
Z.~Wen, D.~Goldfarb, and W.~Yin.
\newblock Alternating direction augmented lagrangian methods for semidefinite
  programming.
\newblock {\em Mathematical Programming Computation}, pages 1--28, 2010.

\bibitem{xu2012alternating}
Y.~Xu, W.~Yin, Z.~Wen, and Y.~Zhang.
\newblock An alternating direction algorithm for matrix completion with
  nonnegative factors.
\newblock {\em Frontiers of Mathematics in China}, 7(2):365--384, 2012.

\bibitem{yang2009fast}
J.~Yang, W.~Yin, Y.~Zhang, and Y.~Wang.
\newblock A fast algorithm for edge-preserving variational multichannel image
  restoration.
\newblock {\em SIAM Journal on Imaging Sciences}, 2(2):569--592, 2009.

\bibitem{yang2009alternating}
J.~Yang and Y.~Zhang.
\newblock Alternating direction algorithms for l1-problems in compressive
  sensing.
\newblock {\em SIAM Journal on Scientific Computing}, 33(1):250--278, 2011.

\bibitem{yang2009efficient}
J.~Yang, Y.~Zhang, and W.~Yin.
\newblock An efficient tvl1 algorithm for deblurring multichannel images
  corrupted by impulsive noise.
\newblock {\em SIAM J. Sci. Comput}, 31(4):2842--2865, 2009.

\bibitem{yang2010fast}
J.~Yang, Y.~Zhang, and W.~Yin.
\newblock A fast alternating direction method for tvl1-l2 signal reconstruction
  from partial fourier data.
\newblock {\em Selected Topics in Signal Processing, IEEE Journal of},
  4(2):288--297, 2010.

\bibitem{yuan2009sparse}
X.~Yuan and J.~Yang.
\newblock Sparse and low-rank matrix decomposition via alternating direction
  methods.
\newblock {\em Preprint}, 2009.

\bibitem{zappella11simul}
L.~Zappella, A.~Del~Bue, X.~Llado, and J.~Salvi.
\newblock Simultaneous motion segmentation and structure from motion.
\newblock In {\em Applications of Computer Vision (WACV), 2011 IEEE Workshop
  on}, pages 679--684, 2011.

\end{thebibliography}

\end{document}